\documentclass[runningheads]{llncs}
\usepackage{ifpdf}
\pdfoutput=1

\usepackage[T1]{fontenc}
\usepackage{nicefrac}

\usepackage[nocompress]{cite}
\usepackage[pdftex]{graphicx}
\usepackage[utf8]{inputenc}
\usepackage[english]{babel}
\usepackage{amsmath}
\usepackage[euler]{textgreek}
\usepackage{dsfont}
\usepackage{afterpage}
\usepackage{booktabs}
\usepackage{multirow}
\usepackage{caption}
\usepackage{subcaption}
\usepackage{float}
\restylefloat{table}
\usepackage{pdflscape}
\usepackage{ragged2e}
\usepackage{amssymb}
\usepackage{subfloat}
\usepackage{natbib}
\usepackage{xcolor}
\usepackage{lipsum}
\usepackage{setspace}
\usepackage{parskip}
\usepackage{fullpage}

\usepackage{fancyhdr}
\usepackage{enumerate}
\usepackage{centernot}
\usepackage{listings}
\usepackage{etoolbox}
\usepackage{mathtools}
\usepackage{algorithm}
\usepackage[noend]{algpseudocode}
\usepackage[toc,page]{appendix}

\renewcommand{\algorithmicrequire}{\textbf{Input:}}
\renewcommand{\algorithmicensure}{\textbf{Output:}}

\DeclareMathOperator*{\argmin}{argmin}

\algrenewcommand\algorithmicrequire{\textbf{Precondition:}}
\algrenewcommand\algorithmicensure{\textbf{Postcondition:}}

\usepackage[colorlinks,allcolors=black]{hyperref} 

\begin{document}

\title{Ridesharing and Fleet Sizing \\ For On-Demand Multimodal Transit Systems}

\author{Ramon Auad\inst{1,2,}\thanks{\texttt{rauad@gatech.edu, ramon.auad@gmail.com, rauad@ucn.cl}, corresponding author} \and
Pascal Van Hentenryck\inst{1,}\thanks{\texttt{pascal.vanhentenryck@isye.gatech.edu}}}
\authorrunning{Auad and Van Hentenryck}

\institute{School of Industrial and Systems Engineering, Georgia Institute of Technology, Atlanta, GA \and Departamento de Ingeniería Industrial, Universidad Católica del Norte, Antofagasta, Chile}
\maketitle              
\begin{abstract}
This paper considers the design of On-Demand Multimodal Transit Systems (ODMTS) that combine fixed bus/rail routes between transit hubs with on-demand shuttles that serve the first/last miles to/from the hubs. The design problem aims at finding a network design for the fixed routes to allow a set of riders to travel from their origins to their destinations, while minimizing the sum of the travel costs, the bus operating costs, and rider travel times. The paper addresses two gaps in existing tools for designing ODMTS. First, it generalizes prior work by including ridesharing in the shuttle rides. Second, it proposes novel fleet-sizing algorithms for determining the number of shuttles needed to meet the performance metrics of the ODMTS design. Both contributions are based on Mixed-Integer Programs (MIP). For the ODMTS design, the MIP reasons about pickup and dropoff routes in order to capture ridesharing, grouping riders who travel to/from the same hub. The fleet-sizing optimization is modeled as a minimum flow problem with covering constraints. The natural formulation leads to a dense graph and computational issues, which is addressed by a reformulation that works on a sparse graph. The methodological contributions are evaluated on a real case study: the public transit system of the broader Ann Arbor and Ypsilanti region in Michigan. The results demonstrate the substantial potential of ridesharing for ODMTS, as costs are reduced by about 26\% with respect to allowing only individual shuttle rides, at the expense of a minimal increase in transit times. Compared to the existing system, the designed ODMTS also cuts down costs by 35\% and reduces transit times by 38\%.

\end{abstract}
\par

\vfill

\pagebreak

\section{Introduction} \label{sec:intro}
\begin{figure}
    \centering
    \includegraphics[width=0.6\textwidth]{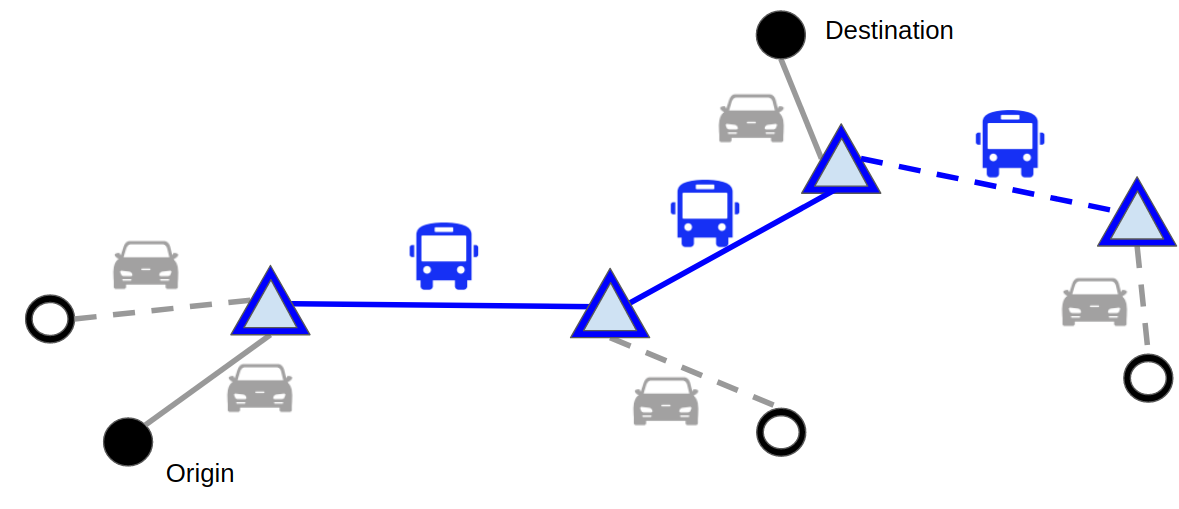}
    \caption{Example of an ODMTS. The Path Followed by a Particular Rider is Represented Using Solid Lines}
    \label{fig:ODMTS_example}
\end{figure}

Recent advances in technology are changing the landscape of city logistics, through the emergence of mobile applications and the concept of shared mobility \citep{kulinska2019development, mccoy2018integrating}. With continuous growth in population and urbanization, city logistics is expected to have a significant societal impact \citep{raghunathan2018seamless, savelsbergh2016, grosse-ophoff2017}. Consequently, it is crucial to envision novel solutions to meet current challenges, and develop cost-effective, environmentally friendly, and socially aware \citep{sampaio2019crowd} transportation systems. One promising solution is the integration of shared mobility and multimodal transportation systems, through a concept known as On-Demand Multimodal Transit Systems (ODMTS) \citep{van2019social}. Figure \ref{fig:ODMTS_example} illustrates the concept of ODMTS, where each passenger travels from an origin to a given destination using routes that combine bus and shuttle legs. This idea presents significant advantages, including improved mobility for those who do not own a vehicle, enhanced first and last mile connectivity, expanded access to public transit systems, and a sustainable business model \citep{kodransky2014connecting, lazarus2018shared, mccoy2018integrating, stiglic2018enhancing, agatz2020make} to name a few. This paper explores this concept by integrating ridesharing into the design of an ODMTS, based on the work by \cite{maheo2017benders}.

The ODMTS Design Problem has been recently studied by \cite{maheo2017benders,dalmeijer2020transfer,basciftci2021capturing, AUAD2021103418}. ODMTS combine fixed-route bus/rail services between selected hubs with on-demand shuttles in charge of transporting riders between the hubs and their origins and destinations. Riders book their trips online (e.g., through a phone app) and are picked up at virtual stops; this facilitates the passenger pickup and avoids delays due to waiting at a customer residence. Riders are also dropped off at a location close to their destination. In between, the trip typically involves one or more bus or rail legs. The tight integration of on-demand shuttle legs with a public transit system may reduce both the travel time and the overall system cost \citep{maheo2017benders,stiglic2018enhancing,van2019social}.

ODMTS also offer opportunities for \emph{ridesharing} and, in particular, grouping, in shared shuttle rides, multiple riders with close origins and heading towards similar destinations. These ride-shared legs further decrease costs and help in improving waiting times during peak times. In particular, they may reduce the number of shuttles, resulting in lower operational costs and, potentially, increased use of public transport \citep{farahani2013review,furuhata2013ridesharing,stiglic2018enhancing}. However, ridesharing in the shuttle legs and fleet sizing have not been considered in the original network design optimization of \citep{maheo2017benders} nor in subsequent work.

This paper aims at addressing this gap: it proposes a two-step framework to (1) integrate ridesharing in the network design optimization and (2) size the shuttle fleet to meet the performance metrics of the proposed design. More precisely, given a set of origin-destination (O-D) pairs and a set of hubs, this paper addresses the problem of designing a cost-efficient ODMTS that considers shared shuttle routes and serves all the transportation requests with the minimum number of shuttles. This relaxes the assumption of \cite{maheo2017benders} that the system always has a shuttle available for serving a request, and provides transit agencies with a precise estimation of the optimal shuttle fleet size, which is critical in practice.

The first step of the framework addresses the ODMTS network design. It selects which pairs of hubs to connect through high-frequency bus routes in order to take advantage of economies of scale, while making use of on-demand shuttles for the first and last legs of the trips. There exists a trade-off between the design of the network and the routing of the passengers via shuttles, since opening a fixed line between hubs requires a significant investment but leads to a low operational cost per trip, whereas shuttle routes have a low up-front cost but a considerably higher cost per ride. To generate shuttle routes to serve riders, this paper employs a pickup and dropoff route enumeration algorithm inspired by \cite{hasan2018community}. The constructed routes are then used as input to a Mixed-Integer Program (MIP) that models the ODMTS design as a Hub-Arc Location Problem (HALP) \citep{campbell2005hub, campbell2005hub2}: the model optimally selects the fixed lines to open and the shuttle routes to serve. The optimal shuttle routes serve as inputs for the second step of the framework, which addresses the fleet-sizing problem. This second optimization model is a minimum cost flow formulation with covering constraints and its coefficient matrix is totally unimodular. It returns the minimum number of shuttles required to complete all the shuttle requests, and the set of requests served by each shuttle. The natural formulation of the fleet-sizing model leads to a dense graph, which raises significant computational difficulties. An improved formulation, that sparsifies the graph, overcomes these limitations. It should be noted that ultimately, the practical goal of solving this problem is to determine a bus network design and a shuttle fleet size; in this regard, the shuttle routing decisions have a direct impact on the ODMTS design and fleet-sizing problem.

The paper presents experimental results to highlight the benefits of incorporating ridesharing into the ODMTS design in terms of operating cost, passengers convenience, and the number of operating shuttles, using real data collected from the Ann Arbor Area Transportation Authority (AAATA). The results demonstrate that ridesharing may improve the ODMTS costs by about 26\% relative to the case where shuttles serve one passenger at a time, while introducing minimal increases in transit times. Compared to the existing system, the findings suggest that the designed ODMTS may also cut down operational daily costs by up to 35\% while potentially reducing transit times by up to 38\%. The paper also validates the model assumptions by performing a sensitivity analysis on key ridesharing parameters, including the impact of estimated arrival times at the last hub and the time window during which riders can be grouped. These results demonstrate the robustness of the proposed framework.

The contributions of this paper can be summarized as follows:
\begin{enumerate}[(i)]

\item it presents a framework to capture ridesharing in the design of an ODMTS, combining a route-enumeration algorithm and a HALP;

\item it formulates the fleet-sizing optimization problem for the on-demand shuttles as a standard vehicle scheduling problem, and proposes an alternative flow formulation that is also totally unimodular but is constructed on a sparse underlying network, significantly enhancing its scalability.

\item it validates the proposed framework through a comprehensive set of experiments using real-world data from a local public transit system, including a sensitivity analysis of the most critical parameters and a comparison with the existing transit system;

\item it presents results that illustrate the potential benefits of ridesharing for ODMTS and the overall benefits in convenience and cost compared to the existing transit system.
  
\end{enumerate} 

The remainder of the paper is organized as follows. Section \ref{sec:lit_rev} covers related literature. Section \ref{sec:net_des} defines the ODMTS network design problem with ridesharing and the underlying assumptions, defines the mathematical notations, presents the route enumeration algorithm, and proposes the MIP model for finding the optimal network design. Section \ref{sec:fleet} introduces the fleet-sizing problem, its assumptions, and its mathematical model. Section \ref{sec:exp_res} reports the numerical results for the case study in Ann Arbor and Ypsilanti, in Michigan, USA. Finally, Section \ref{sec:concl} presents the concluding remarks and future research directions.

\section{Review of Related Literature}
\label{sec:lit_rev}
In the last decades, there have been considerable research on optimizing the design of urban transportation networks. A comprehensive review of this line of research is offered by \cite{farahani2013review}, who compare this problem to the road network design problem \citep{magnanti1984network, yang1998models, xu2016review} and the public transit network design problem \citep{bourbonnais2019transit, demir2016green, cipriani2012transit} in terms of modeling, solution methods, and the most characteristic constraints. The authors further highlight the effect of computational progress, solution methods, and passengers behavior on the evolution of research in the design of urban mobility systems. Similar points are conveyed in the special issue \cite{campbell2019special}.

The problem studied in this paper concerns urban transportation and has its foundations in the work of \cite{campbell2005hub, campbell2005hub2}, which introduced the HALP and several variants. The HALP decision consists of locating a set of arcs between hubs that yields the optimal total flow cost. Among the presented variants, the HAL4 model is most similar to the framework proposed in this paper as it seeks a connected optimal hub network. However, this paper relaxes some of its key assumptions: (1) it allows paths that directly connect, through an on-demand shuttle service, an origin with its corresponding destination; and (ii) it considers solutions with shuttle paths that start or end at a hub node and visit multiple non-hub nodes.

The formulation of the HALP was motivated as an alternative to the hub location problem (HLP) firstly studied by \cite{o1986location}. The HLP is formulated as a MIP, where each decision variable represents whether a particular node is allocated to a particular hub, and path continuity constraints are imposed. The HLP, however, assumes that the hubs form a complete network; such critical assumption is relaxed in the ODMTS design which focuses instead on determining which hub arcs should be opened to obtain economies of scale.  Both problems have a diversity of applications, including the design of large-scale transportation systems, where there are strong opportunities of cost efficiency through consolidation of passengers \citep{campbell2012twenty,lium2009study}. In particular, \cite{campbell2012twenty} address the origins and evolution of the hub location field, and \cite{farahani2013hub, alumur2008network} present an exhaustive survey on hub location literature.

\par This work is closely related to \cite{dalmeijer2020transfer, maheo2017benders}. \cite{maheo2017benders} introduces the ODMTS design problem as part of the {\sc BusPlus} project, seeking to improve the public transportation in the city of Canberra, Australia. By only considering single-passenger shuttle rides, they formulate the design problem as a MIP and identify a special structure in the formulation suitable to employ a Benders decomposition algorithm that generates multiple cuts per iteration \citep{benders1962partitioning}. In addition, the authors propose a pre-processing step that identifies and filters trips that take a direct shuttle trip at optimality, greatly reducing the size of the model, and further accelerate the solution process through the generation of Pareto-optimal Benders cuts \citep{magnanti1981accelerating}. This work is later extended by \citep{dalmeijer2020transfer} who incorporate the selection of the frequency of each opened bus leg and constraints on the number of transfers using a transfer-expanding graph. This allows a Benders decomposition formulation where the sub-problem solves multiple independent shortest path problems. The authors show the effectiveness of this approach using real data from the city of Atlanta, GA.

\par Unfortunately, the incorporation of ridesharing into the problem modeling breaks the special structure mentioned earlier, and consequently neither the Benders cut disaggregation nor the aforementioned pre-processing steps are possible while preserving optimality. If the design optimization is decomposed into a restricted master problem and a sub-problem as in \cite{maheo2017benders}, the sub-problem linear relaxation no longer has extreme integer points, and hence a standard Benders decomposition does not converge to the true optimal solution. Despite such issue, it is still possible to solve instances of reasonable sizes that consider ridesharing in the shuttle legs by only limiting shuttle routes to the ones that satisfy reasonable real-world criteria (e.g., timing constraints). Enforcing these conditions makes it possible to enumerate all the reasonable routes without incurring excessive running times, even for real cases as the one considered in this paper. The route enumeration algorithm in this paper is inspired by the approach in \cite{hasan2018community}, which studies community-based ridesharing. Both algorithms enumerate shared shuttle routes to connect to/from a specific location (e.g., a job center in \cite{hasan2018community} and hubs in the present paper).  As long as the shuttle capacity is not excessively large and only routes of practical interest are considered, the algorithm generates all the routes of potential interest in reasonably short times.

\par Another related line of work include research in last-mile logistics. \cite{raghunathan2018integrated} optimizes total transit time considering the joint schedule of passengers that make use of mass transportation (represented by a train) and fixed shuttle capacity. In this setting, passengers take the train at a given time from a particular station to a common hub, from where they are consequently grouped in shuttle rides that drop them at their final destinations. Under specific conditions, they characterize a set of optimal solutions and further propose a heuristic method that exploits such solution structure. In a later work \cite{raghunathan2018seamless}, the authors propose a more general model that optimizes a convex combination of two objectives, namely the total transit time and the number of shuttle trips. Combining decision diagrams and branch-and-price, they are able to solve real-world instances to optimality in very short times. Additionally, a generalization of this study with uncertainty in the schedule of a subset of passengers can be found in \cite{serra2019last}. The key assumptions present in all these papers are (i) a single hub from where shuttle rides start; (ii) all the stations from where passengers take the initial leg are fixed and visited sequentially; and (iii) any shuttle route has a unique stop and every customer in the ride is dropped off at this stop. This paper relaxes some of these assumptions by considering trips with first and last shuttle legs and middle legs in fixed bus routes, and by making the bus network design a key part of the decision problem. Moreover, shuttle routes may perform multiple intermediate stops to serve riders with different origins and destinations, potentially requiring fewer shuttles to serve all the requests.

\par Fleet size optimization is a fundamental problem for a wide range of applications, including transportation systems planning \citep{baykasouglu2019review}, meal delivery \citep{auad2020demandfleet} and airline scheduling \citep{wang2015optimization}. The fleet-size optimization algorithm proposed in this paper is closely related to the vehicle scheduling literature. A thorough survey of this research line is provided by \cite{bunte2009overview}, which explores multiple problem formulations, variants, and practical extensions. The work by \cite{saha1970algorithm} was the first in solving the single-depot variant of the vehicle scheduling problem. The single-depot variant can be solved in polynomial time and can be formulated as an assignment problem \citep{orloff1976route, kim2012school}, a network flow problem \citep{bokinge1980improved, silva1998vehicle}, or a transportation problem \citep{gavish1979approach, auad2020demandfleet} (unlike its multi-depot counterpart, which is proven to be NP-hard by \cite{bertossi1987some}). \cite{bokinge1980improved} further propose a long arc elimination routine that relies on the value of a depot travel time parameter set by decision makers. This paper proposes an alternative arc elimination algorithm that eliminates transitive arcs from the underlying network and considerably enhances the scalability of the algorithm, while guaranteeing that the optimal solution is not lost. This is especially important in the context of ODMTS design, where the system requires to simultaneously complete a considerable number of shuttle routes.

\section{Network Design with Ridesharing} \label{sec:net_des}

This paper considers an On-Demand Multimodal Transit System (ODMTS) which is composed a set of fixed high-frequency bus lines to serve high density stretches and a set of responsive, on-demand shuttles to serve the first/last miles and act as feeders to/from the fixed routes. The fixed route component addresses congestion and economy of scale, while the on-demand shuttles focus on the first/last mile problem that typically plagues conventional transit systems. In an ODMTS, a trip is typically composed of multiple legs, e.g., a passenger first takes a shuttle to connect to the fixed route network, then traverses the bus network, possibly visiting one or more intermediate stops, and finally takes a second shuttle to connect from the fixed network to the final destination.

The ODMTS design problem considered in this paper consists of designing the bus network that, when integrated with on-demand shuttle routes to serve the first/last miles, minimizes costs and maximizes convenience. In particular, the planner must select which bus lines to open among the candidate arcs, each of which has a fixed opening cost representing the cost of operating high-frequency buses along the corresponding arc during the operating time. The goal is to jointly minimize the total cost of the system, i.e., the fixed cost of operating the bus lines and the variable cost for each shuttle trip, and the inconvenience of the passengers, i.e., the transit time from origin to destination. In addition, the bus network design includes the possibility of passengers sharing shuttle trips, i.e., consolidating multiple passengers in shuttle routes both inbound and outbound to the bus transportation network, up to the shuttle capacity. Shared routes may provide a substantial reduction in the number of shuttles and the total variable cost corresponding to the shuttle rides. The design makes the following assumptions:
\begin{itemize}
    \item Passengers with a common O-D pair and similar departure times are grouped into a single commodity\footnote{The terms \textit{commodity} and \textit{trip} are used interchangeably throughout this work, and refer to a set of passengers with a common O-D pair and similar departure times.} up to the shuttle capacity. If the total number of passengers with a common O-D pair and departure time exceeds the shuttle capacity, the request is split into multiple commodities.
    \item Shuttle routes can be of three types: a direct O-D route, a pickup route, or a dropoff route. A {\em direct route} serves a trip from its origin to its destination and has no ride sharing (except if there are multiple riders in the request). A {\em pickup route} starts at a pickup location, may involve multiple intermediate stops to pick up riders at different locations, and drops all of them off together at a particular hub. A {\em dropoff route} starts at a bus hub with a set of passengers on board, makes a set of sequential stops to drop each of them off, and ends at the destination of the last rider.
    \item Shuttle routes may involve multiple passengers, as long as (i) the individual departure times of the passengers included in the shared route fall in a common predefined time window; and (ii) the total time that each involved passenger spends aboard the shuttle does not exceed a predefined time threshold relative to the duration of the direct route.
    \item The transfer times when connecting between buses are assumed to be fixed and identical throughout the bus network.
    \item Bus lines are only between hubs; a bus that traverses an open line from a hub $h$ to another hub $l$ does not perform any intermediate stops.
\end{itemize}
This work approximates\footnote{Determining the true benefit would require solving a dynamic online variant of the problem where trips information is initially unknown and revealed over time.} the benefits of considering ridesharing at the shuttle legs by:
\begin{enumerate}
\item Solving, for a given set of trips, a static version of the design problem that determines the optimal bus network and associated shuttle routes to be followed by each commodity;
\item Solving the fleet-sizing problem to calculate the minimum number of shuttles required to serve every shuttle leg.
\end{enumerate}
It is important to note that, in this framework, the inclusion of shared shuttle rides serves to guide both the bus network design and the shuttle fleet size, which comprise the practical purpose of solving the studied problem.

\subsection{Problem Formulation} \label{inputs}

The input of the design problem contains the following elements:
\begin{enumerate}[(i)]
    \item a complete graph $G$ with a set $N$ of nodes, where the nodes represent virtual stops and the arcs represent links between them;
    \item a subset $H\subseteq N$ of nodes that are designated as bus hubs;
    \item time and distance matrices $T$ and $D$ that respect the triangle inequality but can be asymmetric: for each $i,j\in N$, $T_{ij}$ and $D_{ij}$ denotes the time and distance from node $i$ to $j$, respectively;
    \item a set $C$ of commodities (trips) to serve: each commodity $r\in C$ is characterized by an origin $or(r)$, a destination $de(r)$, a number of passengers $p(r)$, and a departure time $t_0(r)$;
    \item A time horizon $[T_{min}, T_{max}]$ during which departures occur, i.e., $t_0(r)\in [T_{min}, T_{max}], \forall r\in C$.
\end{enumerate}

\noindent
The ODMTS problem jointly optimizes the fixed cost of opening bus lines, a distance-based cost incurred by the system, and the inconvenience of passengers measured in terms of travel time. The distance cost is computed by multiplying the travel distance by the corresponding shuttle and bus variable costs. To capture costs and inconvenience in a single cost function, the model uses a factor $\alpha$ that balances traveled distance and rider inconvenience, the latter represented as the total travel time incurred by passengers (including waiting times prior to boarding a shuttle and a bus). The objective function is thus the sum of the total inconvenience multiplied by $\alpha$ and the operational cost multiplied by $(1 - \alpha)$. Higher values of $\alpha$ give higher priority to minimizing inconvenience, while lower values translate into an optimal solution that primarily seeks to minimize costs. The following nomenclature is used to compute the total cost:
\begin{itemize}
    \item $K$: the shuttle passenger capacity;
    \item $c$: the variable cost per kilometer of a shuttle;
    \item $b$: the variable cost per kilometer of a bus;
    \item $n$: the number of bus trips for the entire planning horizon across a given opened bus line (assumed to be the same for each line);
    \item $S$: the fixed waiting time incurred by a passenger seeking a bus at a bus hub, from the moment she arrives at the hub until she boards a bus.
\end{itemize}

\noindent
The cost function associated with each mode of transportation accurately captures its characteristics. For buses, let $BL \doteq \{(h,l) \in H\times H : h\neq l\}$ be the set of possible bus lines that can be opened. The decision of opening a bus line $(h, l)$ requires a cost equivalent to the cost of performing $n$ bus trips during a time period of interest from $h$ to $l$ without intermediate stops (thus this cost is modeled as a one-time setup payment). More precisely, for any $(h, l) \in BL$, the cost of opening a bus line from $h$ to $l$ during a time period of interest is explicitly given by
\begin{align*}
    \beta_{hl} \doteq (1 - \alpha) b \cdot n \cdot D_{hl}
\end{align*}
Once bus line $(h, l)$ is opened, the cost incurred by a passenger from using such line is the associated converted inconvenience, i.e.,
\begin{align*}
    \gamma_{hl} \doteq \alpha (T_{hl} + S)
\end{align*}
For a commodity $r\in C$, since waiting and travel times are incurred by each passenger, the inconvenience cost of the $p(r)$ riders using bus line $(h, l)\in BL$ is computed as
\begin{align*}
    \gamma_{hl}^r \doteq p(r)\cdot \gamma_{hl}
\end{align*}
This definition assumes that buses have infinite capacity,
which means that the $p(r)$ riders can always follow the same multi-modal route. 
\par 
The cost of commodity $r\in C$ taking a direct O-D shuttle route is given by
\begin{align*}
    c^{direct}_r \doteq p(r) \cdot \left( (1 - \alpha) c\cdot D_{or(r), de(r)} + \alpha T_{or(r), de(r)}\right)
\end{align*}
On the other hand, every non-direct shuttle route\footnote{Throughout the paper, the term ``shuttle route'' is mainly defined in terms of the sequence of passengers that are picked up or dropped off by the shuttle, in contrast to the standard characterization of routes in terms of the locations that it visits.} $\omega$ is characterized by
\begin{itemize}

\item $k_\omega$: the number of commodities served by route $\omega$.
  
\item $\boldsymbol{r}^\omega$: a vector of commodities $(r_1^\omega, r_2^\omega, \dots, r_{k_\omega}^\omega)$ served by a shuttle following route $\omega$, where $r_j^\omega$ corresponds to the $j$-th commodity picked up (dropped off) in a pickup (dropoff) shuttle route.
  
\item $h_\omega$: the bus hub associated with route $\omega$; in pickup routes, $h_\omega$ corresponds to the route ending point; in dropoff routes, $h_\omega$ corresponds to the starting point of the route; direct O-D routes do not involve hubs and so this parameter does not apply.

\item $\boldsymbol{\xi}^\omega$: a time vector $(\xi^\omega_1, \xi^\omega_2, \dots, \xi^\omega_{k_\omega})$ where $\xi_j^\omega$ denotes the total time that commodity $r_j^\omega$ incurs to complete route $\omega$. For pickup routes, $\xi_j^\omega$ corresponds to the time from departure time $t_0(r_j^\omega)$ to when $r_j^\omega$ leaves the shuttle. For dropoff routes, and assuming the route $\omega$ starts at a hub $h$, $\xi_j^\omega$ represents the time period between the arrival of commodity $r_j^\omega$ to hub $h$ (possibly having to wait for the arrival of more commodities to location $h$ prior to starting the route), and the time at which $r_j^\omega$ is dropped off by the shuttle at its final destination (the computation of $\boldsymbol{\xi}^\omega$ for dropoff routes is discussed in more detail in Section \ref{sec:route_enum}).
    
\item $p_\omega$: the total number of passengers picked up (dropped off) by a shuttle following route $\omega$, with $$p_\omega \doteq \sum_{j=1}^{k_\omega} p(r_j^\omega)$$
    
\item $\mathcal{A}_\omega$:  the set of arcs $(i, j)\in N\times N$ traversed by shuttle route $\omega$.

\item $d_\omega$: the total distance driven by a shuttle following route $\omega$, i.e., $$d_\omega \doteq \sum_{(i,j)\in \mathcal{A}_\omega}D_{ij}$$
    
\item $c_\omega$: the total cost (combining distance cost and inconvenience) incurred by a shuttle following route $\omega$, computed as $$c_\omega \doteq (1 - \alpha) c \cdot d_\omega + \alpha \sum_{j=1}^{k_\omega}p(r_j^\omega) \cdot \xi_j^\omega$$

\end{itemize}

\noindent
Direct shuttle routes result in a lower inconvenience, but routes serving multiple trips have lower costs.

\subsection{The MIP Model} \label{sec:mip_model}

This section presents the MIP model associated with the design of the ODMTS. The MIP model receives as input a set of shuttle routes and uses the following notations:
\begin{itemize}
\item $\Omega^-_r$: the set of pickup routes $\omega$ such that $r \in \boldsymbol{r}^\omega$ for commodity $r\in C$;
\item $\Omega^+_r$: the set of dropoff routes $\omega$ such that $r \in\boldsymbol{r}^\omega$ for commodity $r\in C$.
\end{itemize}
The set of pickup routes is denoted by $\Omega^- \doteq \bigcup_{r\in C} \Omega_r^-$ and the set of dropoff routes by $\Omega^+ \doteq \bigcup_{r\in C} \Omega_r^+$. The construction of these routes is discussed in Section \ref{sec:route_enum}.

The MIP model considers two interacting decisions: it determines (i) which bus lines to open, and (ii) which route riders follow from their origin to their destination, either using a direct route or multi-modal routes combining shuttle and bus legs. Multi-modal routes can only use opened bus legs. The MIP formulation models these decisions using the following binary decision variables: $z_{h, l} = 1$ iff bus line $(h, l) \in BL$ is selected to be opened; $y_{h, l}^r = 1$ iff riders in $r\in C$ take bus line $(h, l)\in BL$; $x_{\omega} = 1$ iff shuttle route $\omega\in \Omega^{-}\cup \Omega^{+}$ is selected to be served; $\eta_r = 1$ iff riders in $r\in C$ take a direct shuttle route from $or(r)$ to $de(r)$.

Model \eqref{mip_model} presents the MIP model. Objective \eqref{tot_cost} minimizes the total cost, which includes the routing costs (the cost and inconvenience of direct and multi-modal routes) and the cost of opening bus lines.  Constraints \eqref{topology} enforce a weak connectivity on the resulting bus network, requiring that, for each hub $h\in H$, the number of opened bus lines inbound to $h$ must match the number of outbound opened lines. As mentioned in \cite{maheo2017benders}, although \eqref{topology} by itself does not theoretically guarantee full connectivity of the resulting bus network, in practice, the spatial distribution of the origins and destinations makes this set of constraints sufficient for this purpose.  Constraint sets \eqref{cover_p} and \eqref{cover_d} guarantee that each commodity $r\in C$ is both picked up at its origin and dropped off at its destination, either by a direct or a shared route.  Constraints \eqref{bus_leg_available} restrict bus legs to only use opened bus lines, and Constraints \eqref{hub_flow_conservation} enforce the flow conservation constraints at each hub.

\begin{subequations}\label{mip_model}
\begin{align}
	\min \;\;& \sum_{(h, l) \in BL} \beta_{h, l} z_{h, l} + \sum_{r \in C} \bigg(c^{direct}_r \eta_r + \sum_{\omega\in \Omega^-_r}c_{\omega} x_{\omega} + \sum_{\omega\in \Omega^+_r}c_{\omega} x_{\omega} + \sum_{(h, l) \in BL} \gamma^r_{h, l} y_{h, l}^r\bigg) \label{tot_cost}
	\\
	\text{s.t.} \;\; & \sum_{l \in H} z_{h, l} = \sum_{l\in H} z_{l, h} \qquad \forall h\in H 
	\label{topology} 
	\\
	& \eta_r + \sum_{\omega\in\Omega^-_r} x_{\omega} \geq 1 \qquad \forall r\in C 
	\label{cover_p}
    \\
    & \eta_r + \sum_{\omega\in\Omega^+_r} x_{\omega} \geq 1 \qquad \forall  r\in C 
    \label{cover_d}
    \\
    & y^r_{h, l} \leq z_{h, l} \qquad \forall (h, l) \in BL, \quad \forall r\in C 
    \label{bus_leg_available}
    \\
    & \sum_{l\in H} y^r_{l,h} + \sum_{\substack{{\omega}\in \Omega^-_r \\ \text{if } h_{\omega} = h}} x_{\omega} = \sum_{l\in H} y^r_{h,l} + \sum_{\substack{\omega\in \Omega^+_r \\ \text{if } h_{\omega} = h}} x_{\omega} \qquad\quad \forall r \in C, \quad \forall h \in H  
    \label{hub_flow_conservation}
    \\
    & z_{h, l}, y_{h, l}^r, x_{\omega}, \eta_r \in \{0,1\}
\end{align}
\end{subequations}

\subsection{The Route Enumeration Algorithm}
\label{sec:route_enum}

This section describes the generation of the shared routes used as
inputs for the Model \eqref{mip_model}.

\subsubsection{Practical Considerations.}

The algorithm restricts attention to routes of practical interest, using a route duration threshold $\delta > 0$, a consolidation time bucket length $W > 0$, a set of feasible first hubs $H_r^-\subseteq H$ for trip $r\in C$ to enter the bus network, and a set of feasible last hubs $H_r^+$ for trip $r\in C$ to exit the bus network. Consider a sequence of $m\geq 1$ commodities $(r_1,r_2,\dots,r_{m})$ and a hub $h\in H$. In order for the route enumeration algorithm to define a route $\omega$ with $\boldsymbol{r}^\omega = (r_1, r_2,\dots, r_m)$ and $h_\omega = h$, $\omega$ must satisfy three conditions:

\begin{enumerate}
    \item if $\omega$ is a pickup route, then $h\in H_{r_j^{\omega}}^-$ and $\xi_j^{\omega} \leq (1 + \delta) \cdot  T_{or(r^{\omega}_j), h}$, $j\in\{1,2,\dots,m\}$;
    \item if $\omega$ is a dropoff route, then $h\in H_{r_j^{\omega}}^+$ and $\xi_j^{\omega} \leq (1 + \delta) \cdot T_{h, de(r^{\omega}_j)}$, $j\in\{1,2,\dots,m\}$;
    \item $p_\omega \leq K$.
\end{enumerate}

\noindent
Condition 1 requires that hub $h$ is a feasible first hub for all trips in the route, i.e. $h\in H_r^-, \forall r\in \boldsymbol{r}^\omega$, and that the total time spent by commodity $r_j$ in a shared pickup route towards hub $h$ does not exceed $(1 + \delta)$ times the duration of the direct shuttle route from $or(r_j)$ to $h$; and condition 2 imposes similar requirements for dropoff routes. Condition 3 enforces that the number of riders served by a route cannot exceed the shuttle capacity $K$.

Ride-shared routes should only consider riders with close departure times. The operating time horizon $[T_{min}, T_{max}]$ is partitioned into $\lceil \tfrac{T_{max} - T_{min}}{W} \rceil$ time buckets of $W$ minutes. A set of commodities can be served by a shuttle route only if their departure times lie in one of these $W$-minute time buckets. Pickup routes can easily be consolidated based on the departure times of their riders (i.e., $t_0(r)$, $r\in C$). However, dropoff routes raise an interesting issue since the arrival of riders at their starting hubs requires an ODMTS design. To overcome this difficulty, for each commodity $r\in C$ and each hub $l$, the algorithm approximates the time $t_1(r, l)$ when the $p(r)$ riders may reach hub $l$ in their path toward their final destination $de(r)$; this approximation is then used to decide which commodities can be grouped together in a dropoff route. This estimation is computed as the average of the total travel times obtained from each of the $|H_r^-|$ paths that start at $or(r)$ at time $t_0(r)$, travel by shuttle to one of the $|H_r^-|$ existing feasible first hubs, and then take a bus leg to $l$, i.e.,
\begin{align*}
    t_1(r, l) \doteq t_0(r) + \tfrac{1}{|H_{r}^-|} \sum_{h\in H_{r}^-} (T_{or(r), h} + S + T_{h, l}).
\end{align*}
Note that the only purpose of this approximation is to decide which riders may be grouped together to avoid the generation of impractical shared routes. As a result, a shuttle route $\omega$ shared by any two commodities $r, s\in C$ must satisfy one of the following timing conditions:
\begin{enumerate}
\setcounter{enumi}{3}
\item if $\omega$ is a pickup route, then there exists $q\in \mathbb{Z}_+$ such that $t_0(r), t_0(s) \in [T_{min} + qW, \min\{T_{min} + (q+1)W, T_{max}\}]$;
\item if $\omega$ is a dropoff route, then there exists $q\in \mathbb{Z}_+$ such that $t_1(r, h_\omega), t_1(s, h_\omega) \in [T_{min} + qW, \min\{T_{min} + (q+1)W, T_{max}\}]$.
\end{enumerate}

\noindent
These considerations are motivated by the fact that riders may not agree to share a shuttle if the shared route results in considerably longer travel or waiting times.

\subsubsection{The Algorithm.}

\begin{algorithm}[!t]
\caption{Pickup Route Enumeration}\label{alg1}
\begin{algorithmic}[1]
\Require Set of commodities $C$, shuttle capacity $K$, sets of feasible first hubs $\{H_{r}^-\}_{r\in C}$, travel time threshold $\delta$
\Ensure For each $r\in C$, set of pickup routes $\Omega^-_r$
\State $\overline{C}_K \gets \emptyset$
\For{$r_1\in C$}
    \State $\Omega_{r_1}^- \gets \emptyset$
    \For{$h\in H_{r_1}^-$} \label{hubs_-}
        \State $\omega_{r_1} \gets \text{ individual pickup route with $\boldsymbol{r}^\omega = (r_1)$ and $h_\omega = h$}$
        \State $\Omega_{r_1}^-\gets \Omega_{r_1}^- \cup \{\omega_{r_1}\}$ \label{ind_routes_-}
        \For{$k\in \{2,\dots, K\}$} \label{cap_-}
            \State $\sigma_{\text{perm}} \gets \{$All $(k-1)$-element permutations of trips $(r_2, \dots, r_k)\in (C\setminus\{r_1\})^k$ such that a route $\omega$ with $\boldsymbol{r}^\omega = (r_1, r_2,\dots, r_k)$ and $h_\omega = h$ satisfies practical conditions 1, 3, and 4\} \label{considering}
            \For {$(r_2, \dots, r_k) \in \sigma_{\text{perm}}$}
                \State $\overline{C}_K \gets \overline{C}_K \cup \{(h, \{r_1,r_2,\dots, r_k\})\}$ \label{line:C_bar_update}
                \State $\omega \gets \text{ pickup route with $\boldsymbol{r}^\omega = (r_1, r_2,\dots, r_k)$ and $h_\omega = h$}$ \label{construct_route}
                \For{$j\in\{1,2,\dots,k\}$}
                    \State $\Omega_{r_j}^{temp} \gets \Omega_{r_j}^{temp} \cup \{\omega\}$\label{include_-}
\EndFor
\EndFor
\EndFor
\EndFor
\EndFor
\For {$(h, \overline{C}) \in \overline{C}_K$} \label{line:pre_filter_1}
    \State $\omega^* = \argmin\{c_{\omega} : \omega\in\bigcap_{r\in\overline{C}} \Omega_r^{temp} \text{ and } h_\omega = h\}$
    \For {$r\in \overline{C}$}
        \State $\Omega_r^- \gets \Omega_r^- \cup \{\omega^*\}$ \label{line:pre_filter_2}
    \EndFor
\EndFor
\end{algorithmic}
\end{algorithm}

This section describes the algorithm to construct the sets of routes $\Omega^-_{r}$ and $\Omega^+_{r}$ for every commodity $r\in C$, considering homogeneous shuttles with fixed capacity $K$. Algorithm \ref{alg1} sketches the enumeration process for $\Omega_r^-$. For each $r_1\in C$ and feasible first hub $h\in$ $H_{r_1}^-$, the algorithm first generates the individual pickup route that travels from $or(r_1)$ to $h$ (lines \ref{hubs_-} - \ref{ind_routes_-}). Then for the multi-passenger routes, it sets commodity $r_1$ as the first pickup in the route and iterates over all the possible permutations of sizes 1 up to $K-1$ of the remaining commodities in $C$, considering only permutations of commodities whose travel time in $\omega$ satisfies Conditions 1, 3, and 4 (\ref{cap_-} - \ref{considering}). For each such permutation, line \ref{line:C_bar_update} stores the pair ($h, \{r_1,\dots, r_k\})$ in the set $\overline{C}$ to later perform a route pre-filtering, and then the algorithm constructs a route $\omega$ that picks up commodities $r_1,r_2,\dots,r_k$ in that order and drops them off at hub $h$ (line \ref{construct_route}) and adds this route to the set of pickup routes $\Omega_{r_j}^-$ of each picked up commodity $r_j, j\in\{1,\dots,k\}$ (line \ref{include_-}). The procedure is repeated by fixing every commodity $r\in C$ to be the first pickup in a route. Note that the enumeration algorithm evaluates $O\left(\tfrac{(|C| - 1)!}{(|C| - K)!}\sum_{r\in C}|H_r^-|\right)$ routes: for each commodity $r \in C$ that is first picked up in a route, the potential shared routes consist of picking up up to $K-1$ of the remaining $|C|-1$ trips in every possible order, and each of these shared routes may end in any of the $|H_r^-|$ feasible first hubs. However, in practice, it is possible to greatly speed up the enumerating process by using a depth-first approach that prunes the search space by exploiting practical conditions 1, 3, and 4.  Additionally, although the algorithm may generate multiple routes that transport the same subset of commodities $\overline{C}$ to a hub $h$ in different pickup orders, only the least cost route among them is of practical interest and selected by the optimization model. Hence, lines \ref{line:pre_filter_1} - \ref{line:pre_filter_2} only keeps, for each set of commodities $\overline{C}$ served together, the least-cost route serving them together into $\Omega^-_r, \;\forall r\in \overline{C}$. This allows to significantly reduce the number of generated routes. The algorithm to construct the sets of dropoff routes $\Omega^+_r$ follows an almost identical sequence of steps as Algorithm \ref{alg1}.

\section{Fleet-Sizing Optimization}
\label{sec:fleet}

This section discusses the fleet-sizing optimization that minimizes the number of shuttles needed by the ODMTS. It starts with a general formulation which is then improved for computational efficiency.

\subsection{General Formulation} \label{sec:fleet_size_1}

Given the set of optimal shuttle routes defined by solution vectors $\eta^*$ and $x^*$ from solving Model \eqref{mip_model}, this section presents a MIP model which minimizes the number of shuttles required to serve all these routes on time.  The input for the fleet-sizing optimization is a set of shuttle routes $\Omega$ obtained by solving the ODMTS design model, i.e., $\Omega = \{\omega \in \Omega^- \cup \Omega^+ : x_\omega^* = 1\}$, where each route $\omega\in \Omega$ is characterized by a start location $\ell_0^\omega$, an end location $\ell_f^\omega$, a start time $\tau_\omega$, and a duration $\Delta_\omega$. In particular,
\begin{itemize}

\item If $\omega$ is a pickup route, then $(\ell_0^\omega, \ell_f^\omega, \tau_\omega, \Delta_\omega) = (or(r^\omega_1), h_\omega,t_0(r_1^\omega),\xi_1^\omega)$. The route starts at location $or(r_1^\omega)$ where the first commodity $r_1^\omega$ is picked up at departure time $t_0(r_1^\omega)$. Moreover, the route ends at the hub $h_\omega$, where all the pickups are dropped off for a total duration of $\xi_1^\omega$.
    
\item If $\omega$ is a dropoff route then $(\ell_0^\omega, \ell_f^\omega, \tau_\omega, \Delta_\omega) = (h_\omega, de(r^\omega_{k_\omega}), \max_j\{t_1(r^\omega_j, h_\omega)\},\xi_{k_\omega} ^\omega)$. The route starts at hub $h_\omega$ when all commodities $\boldsymbol{r}^\omega$ arrive to $h_\omega$ and the start time is computed as $\max_j\{t_1(r_j^\omega, h_\omega)\}$. The route ends at $de(r_{k_\omega}^\omega)$ when the last commodity $r_{k_\omega}^\omega$ is dropped off and hence its duration is $\xi_{k_\omega}^\omega$.

\item Direct O-D routes are also considered in the set of routes $\Omega$. Particularly, for each $r\in C$ such that $\eta^*_r = 1$, $\Omega$ considers $p(r)$ additional individual routes, each with $(\ell_0^\omega, \ell_f^\omega, \tau_\omega, \Delta_\omega) = (or(r), de(r), t_0(r),T_{or(r), de(r)})$. Each of these routes starts at location $or(r)$ at $t_0(r)$ and travels directly to $de(r)$, arriving at time $t_0(r) + T_{or(r), de(r)}$.
\end{itemize}

The fleet-sizing optimization first builds a directed graph $\mathcal{G} = (\mathcal{V}, \mathcal{A})$ with a unique source node $s$ and sink node $s'$, and where each node in $\mathcal{V}\setminus\{s,s'\}$ uniquely represents a shuttle route in $\Omega$. As a result, the presentation in this section uses ``node'' and ``route'', as well as $\mathcal{V}\setminus\{s,s'\}$ and $\Omega$, interchangeably. The source connects to every node $\omega\in\Omega$ through an arc $(s,\omega)$, and every node $\omega\in\Omega$ connects to the sink $s'$ via an arc $(\omega, s')$. Furthermore, for each pair of routes $\omega, \mu\in\mathcal{V}\setminus\{s,s'\}$, there is an arc $(\omega, \mu) \in \mathcal{A}$ when a single shuttle may feasibly serve routes $\omega$ and $\mu$ in that order, i.e., when $\tau_\omega + \Delta_\omega + T_{\ell_f^\omega, \ell_0^\mu} \leq \tau_\mu$. The construction algorithm for $\mathcal{G}$ is shown in Algorithm \ref{fleet_alg}.

\begin{algorithm}[!t]
\caption{The Fleet-Sizing Graph Construction.} \label{fleet_alg}
\begin{algorithmic}[1]
\Require Set of routes $\Omega$.
\Ensure Task network $\mathcal{G} = (\mathcal{V}, \mathcal{A})$.
\State Let $s$ and $s'$ be the source and sink nodes, respectively.
\State $\mathcal{V}\gets \Omega \cup \{s, s'\}$, $\mathcal{A} \gets \emptyset$ 
\For{$\omega\in \Omega$}
\State $\mathcal{A} \gets \mathcal{A} \cup \{(s, \omega), (\omega, s')\}$
\For{$\mu \in \{\omega' \in \Omega : \tau_{\omega'} > \tau_\omega\}$}
\If{$\tau_\omega + \Delta_\omega + T_{\ell_f^\omega, \ell_0^\mu} \leq \tau_\mu$}
\State $\mathcal{A} \gets \mathcal{A} \cup \{(\omega, \mu)\}$
\EndIf
\EndFor
\EndFor
\end{algorithmic}
\end{algorithm}

The fleet-sizing optimization uses a binary decision variable $v_{\omega,\mu}$ for each route pair $(\omega, \mu)\in \mathcal{A}$ whose value is 1 iff a shuttle serves route $\mu$ immediately after serving route $\omega$. Let $\delta^-_\omega \doteq \{\omega' \in \mathcal{V} : (\omega', \omega) \in \mathcal{A}\}$, and $\delta^+_\omega \doteq \{\omega' \in \mathcal{V} : (\omega, \omega') \in \mathcal{A}\}$. Model \eqref{mip:fleet_size_base} presents the MIP model to minimize the fleet size needed to serve all the selected shuttle routes associated with a given hub. Objective \eqref{fleet_size} captures the number of shuttles needed to cover all routes as the total flow from the source $s$ to any other node. Constraints \eqref{any_copy} require that every route $\omega\in\Omega$ is visited by one unit of flow, and Constraints \eqref{flow_conserv} enforce flow conservation at all nodes other than the source and sink. The coefficient matrix of the fleet-sizing model is totally unimodular: since the right-hand side is integer, the model can be formulated as a linear program.

\begin{subequations} \label{mip:fleet_size_base}
\begin{align}
    \min \;\; & \sum_{\omega\in \delta^+_s} v_{s,\omega} \label{fleet_size}
    \\
    s.t. \;\; &\sum_{\mu\in \delta^-_\omega} v_{\mu, \omega} = 1, \qquad \forall \omega\in \Omega
    \label{any_copy}
    \\
    & \sum_{\mu\in \delta^-_\omega} v_{\mu,\omega} = \sum_{\mu\in \delta^+_\omega}v_{\omega,\mu}, \qquad \forall \omega\in\Omega
    \label{flow_conserv}
    \\
    & v_{\omega,\mu} \in\{0,1\}, \quad \forall (\omega, \mu) \in \mathcal{A}
\end{align}
\end{subequations}

\begin{figure}[!t]
    \centering
    \includegraphics[width=\textwidth]{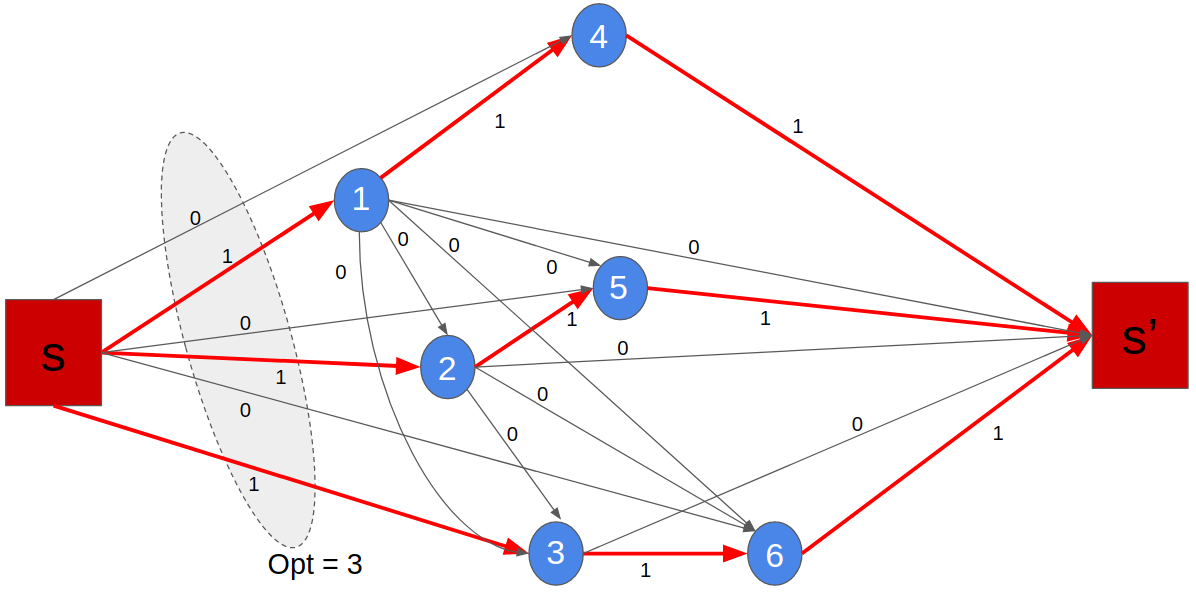}
    \caption{A Fleet-Sizing Graph with $|\Omega| = 6$ Routes and an Optimal Fleet Size of 3.}
    \label{fig:fleet_network}
\end{figure}

Figure \ref{fig:fleet_network} provides an example of the output of Algorithm \ref{fleet_alg} for a simple instance with $\Omega = \{1,2,\dots,6\}$. Routes are indexed based on the start time, with smaller indices implying earlier start times.\footnote{Note that arc $(2,4)$ is not defined even though $(2, 3)$ is and $\tau_3 < \tau_4$. This is because the repositioning time $T_{\ell_f^2,\ell_0^4}$ is too long for a single shuttle to serve routes 2 and 4 in sequence.} Given that the arcs represent all the feasible sequential completions of routes, at least 3 shuttles are required to complete all the requests on time. The solution of this formulation also specifies the sequence of routes each shuttle serves, which opens the possibility to optimize other objectives that depend on this information (e.g., driven distance, total travel time).

\subsection{A Sparse Fleet-Sizing Formulation}

Experimental results on solving the fleet-sizing MIP Model \eqref{mip:fleet_size_base} indicated that practical case studies create an excessive number of feasible arcs, slowing down the solving process considerably due to the large number of variables and significant memory consumption. To overcome these computational issues, this section introduces an arc-filtering procedure that results in a significantly sparser fleet-sizing graph. The key idea underlying the filtering is the fact that shuttles are formulated as a flow and that computing the optimal fleet size only requires ensuring that every node is visited by at least one shuttle; as long as this is satisfied, it is not necessary to explicitly define all the arcs between routes. As a result, it is possible to reduce the number of arcs by removing transitive arcs between routes: if arcs $(\omega_1, \omega_2)$ and $(\omega_2, \omega_3)$ are defined, then it is not necessary to define arc $(\omega_1, \omega_3)$ even though it represents a feasible service sequence. Instead, it is sufficient to remove the capacity limit of arcs $(\omega_1, \omega_2)$ and $(\omega_2, \omega_3)$ and to allow ``multiple shuttles'' to traverse the arcs. Once the new formulation is solved, it is possible to recover the path followed by each shuttle. At termination, the algorithm produces a set of routes whose cardinality matches to the optimal fleet size.

\begin{algorithm}[!t]
  \caption{The Sparse Fleet-Sizing Graph Construction.}
  \label{fleet_alg_tuned}
\begin{algorithmic}[1]
\Require List of routes $\Omega$ sorted by earliest start route.
\Ensure Task network $\mathcal{G} = (\mathcal{V}, \mathcal{A})$.
\State Let $s$ and $s'$ be source and sink nodes, respectively.
\State $\mathcal{V}\gets \Omega \cup \{s, s'\}$, $\mathcal{A} \gets \emptyset$ 
\For{$\omega\in \Omega$}
\For{$\mu \in \Omega_\omega$}
\If{$\{\omega' \in \Omega_{\omega}: \mu \in \Omega_{\omega'}\} = \emptyset$} \label{filter_1}
\State $\mathcal{A} \gets \mathcal{A} \cup \{(\omega, \mu)\}$ \label{filter_2}
\EndIf
\EndFor
\If{$\delta^-_\omega= \emptyset$}
\State $\mathcal{A} \gets \mathcal{A} \cup \{(s, \omega)\}$
\EndIf
\EndFor
\For{$\omega\in \Omega$}
\If{$\delta_\omega^+ = \emptyset$}
\State $\mathcal{A} \gets \mathcal{A} \cup \{(\omega, s')\}$
\EndIf
\EndFor
\end{algorithmic}
\end{algorithm}

\par To formulate the new graph construction algorithm, consider each route $w\in \Omega$ and let $\Omega_{\omega} \doteq \{w'\in\Omega : \tau_\omega + \Delta_\omega + T_{\ell_f^\omega, \ell_0^{\omega'}} \leq \tau_{\omega'}\}$ be the set of routes that may be served immediately after $\omega$ with the same shuttle. The modified network construction procedure is then presented in Algorithm \ref{fleet_alg_tuned}. For routes $\omega\in \Omega$ and $\mu\in \Omega_\omega$, the arc $(\omega, \mu)$ is created only if no intermediate route $\omega'$ exists such that $\omega'\in \Omega_\omega$ and $\mu \in \Omega_{\omega'}$, as stated in lines \ref{filter_1} and \ref{filter_2}.

\begin{subequations} \label{mip:fleet_size_tuned}
\begin{align}
    \min \;\; & \sum_{\omega\in \delta^+_s} v_{s,\omega} \label{fleet_size_2}
    \\
    s.t. \;\; & \sum_{\mu\in \delta^-_\omega} v_{\mu, \omega} \geq 1, \qquad \forall \omega\in\Omega
    \label{any_copy_2}
    \\
    & \sum_{\mu\in \delta^-_\omega} v_{\mu,\omega} = \sum_{\mu\in \delta^+_\omega}v_{\omega,\mu}, \qquad \forall \omega\in\Omega
    \label{flow_conserv_2}
    \\
    & v_{\omega,\mu} \in \mathbb{Z}_+, \quad \forall (\omega, \mu) \in \mathcal{A} \label{domain_2}
\end{align}
\end{subequations}

Given this new fleet-sizing graph, it is possible to define a new optimization model with the following decision variables: variable $v_{\omega, \mu}$ represents the number of shuttles traversing arc $(\omega, \mu)$. Model \eqref{mip:fleet_size_tuned} presents the sparse fleet-sizing optimization model. Objective \eqref{fleet_size_2} minimizes the total number of shuttles used to complete all the routes. Constraints \eqref{any_copy_2} ensure that every node is visited by at least one shuttle; this is a relaxation with respect to Model \eqref{mip:fleet_size_base} that is necessary due to the more limited number of arcs in the sparser graph. Constraints \eqref{flow_conserv_2} enforce flow conservation at all nodes, and Constraints \eqref{domain_2} admit uncapacitated flows but requires them to take integer values.

\begin{figure}[!t]
    \centering
\includegraphics[width=\textwidth]{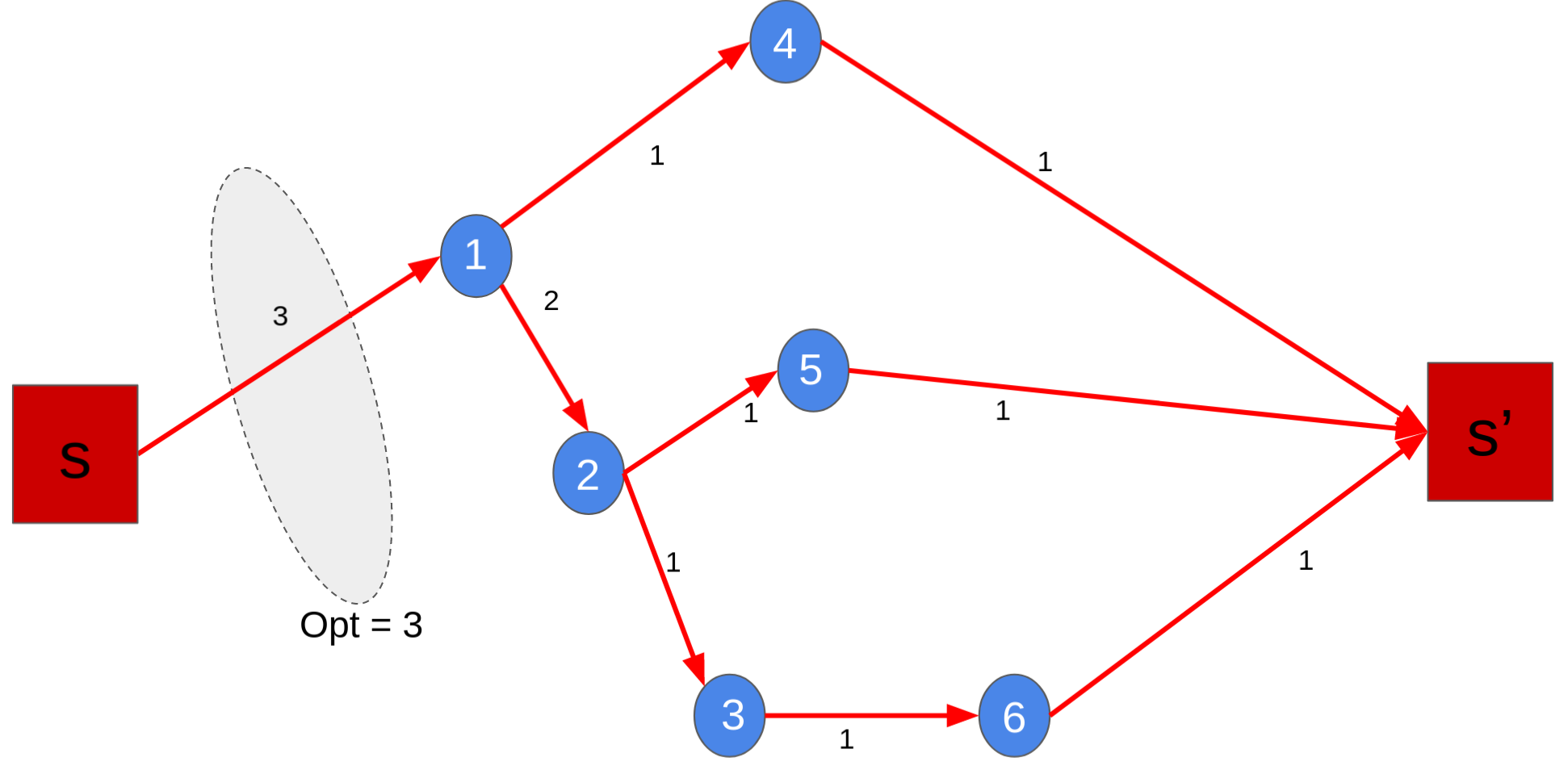}
    \caption{Sparse fleet-sizing graph for the example with $|\Omega| = 6$ shuttle routes. The new formulation correctly concludes that the optimal fleet size is 3, with set of shuttle schedules $P=\{(1,4), (2,5),(3, 6)\}$.}
    \label{fig:fleet_network_2}
\end{figure}

The sparse graph construction for the example from Section \ref{sec:fleet_size_1} is illustrated in Figure \ref{fig:fleet_network_2}. Despite the notorious simplicity of the new network compared to the one from Section \ref{sec:fleet_size_1}, the new underlying optimization model is still able to determine the correct optimal fleet size. The resulting set of shuttle schedules can be either $P=\{(1,4), (2,5), (3,6)\}$ or $P=\{(1,4), (2,3,6), (5)\}$. An algorithm to compute the schedule to be followed by each shuttle is provided in Appendix \ref{appA}.

\section{Experimental Results} \label{sec:exp_res}

This section describes computational results for an ODMTS with ridesharing to illustrate the technical results of the paper. It reports results obtained by solving the ODMTS design Model \eqref{mip_model} and the sparse fleet-sizing Model \eqref{mip:fleet_size_tuned} using data from a real case study concerning the AAATA public transit system for the broad region of Ann Arbor and Ypsilanti in Michigan. The considered transit system network comprises a set $N$ of $1267$ virtual stops, with a subset $H$ of 10 bus hub candidates for the final ODMTS bus network design. With the exception of Section \ref{sec:exp:transit_system_eval}, results are reported for a base instance that considers historical data of 6606 riders who used the public transit system between 6:00 am and 10:00 am of a particular weekday (respectively denoted as $T_{min} = 0$ and $T_{max} = 240$ minutes). The values for the passenger parameters are summarized in Table \ref{tab:route_params}. In general, the grouping of riders into routes assumes a time bucket of length $W = 3$ minutes. For instance, riders who request a ride between 6:00:00 AM and 6:02:59 AM are in principle eligible for sharing a shuttle route, and so on for each subsequent $3$-minute bucket (Section \ref{sec:exp:W} explores the effect of varying this value). The grouping of passengers also assumes a shuttle route duration threshold with $\delta = 50\%$, and the sensitivity of the system to variations in this parameter is explored in Section \ref{sec:delta_sensitivity}. Additionally, to prevent excessively long shuttle routes, the analysis assumes that, for each trip $r$, the set of feasible first hubs $H_r^-$ is given by the 3 hubs closest to $or(r)$, and the set of feasible last hubs $H_r^+$ consists of the 3 hubs closest to $de(r)$. The shuttle capacity used in the experiments varies per section: Sections \ref{sec:exp:base_case} and \ref{sec:exp:W} explore values $K\in\{1,2,3,4\}$ and analyze their effect on the system performance, while Sections \ref{sec:exp:t0+} and \ref{sec:exp:transit_system_eval} consider a unique capacity value of $K=3$.

Table \ref{tab:cost_params} shows the cost-related parameters. The cost structure considers a shuttle variable cost per kilometer of $c = \$1.00$ and a bus variable cost of $b = \$3.75$ per kilometer. These costs assume that (1) shuttles drive at an average speed of 17 miles per hour and cost $\$27.00$ per hour; and (2) buses drive at an average speed of 12 miles per hour and cost $\$72.00$ per hour. Furthermore, the frequency in any opened bus line is set to 4 buses per hour, which translates into an average bus transfer time of $S = 7.5$ minutes and a total of $n = 16$ buses per opened line in the operating period between 6:00am and 10:00am. Traveled distance and riders' inconvenience are balanced by using a value of $\alpha = 10^{-3}$.

The results are presented in the form of six key metrics: the total operating cost of the system in dollars, the average inconvenience of the riders in minutes, the optimal network design, the average shuttle utilization as the number of riders per shuttle route, the number of riders who use direct O-D routes, and the optimal fleet size required to serve all the requests.

\begin{table}[!t]
\caption{Routing-Related Parameter Configurations for Experiments.}
\label{tab:route_params}
\centering
\begin{tabular}{cccccc}
\hline
Section & $K$ & $\delta$ & $W$ (minutes) & $|H_r^-|$ & $|H_r^+|$ \\ \hline
\ref{sec:exp:base_case} & $\{1,2,3,4\}$ & 50\% & 3 & 3 & 3\\
\ref{sec:exp:W} & $\{1,2,3,4\}$ & 50\% & $\{1,3,5\}$ & 3 & 3\\
\ref{sec:delta_sensitivity} & 3 & \{50\%, 100\%, 150\%\} & 3 & 3 & 3\\
\ref{sec:exp:t0+} & 3 & 50\% & 3 & 3 & 3\\
\ref{sec:exp:transit_system_eval} & 3 & 50\% & 3 & 3 & 3\\ \hline
\end{tabular}
\end{table}

\begin{table}[!t]
\caption{Cost-Related Parameter Configuration for Experiments.}
\label{tab:cost_params}
\centering
\begin{tabular}{ccccc}
\hline
$S$ (min) & $\alpha$ & $b$ & $c$ & $n$ \\ \hline
7.5 & $10^{-3}$ & \$3.75 & \$1.00 & 16 \\
\hline
\end{tabular}
\end{table}

\subsection{The Case Study} \label{sec:exp:base_case} 

\begin{figure}[!t]
    \centering
    \includegraphics[width=0.7\textwidth]{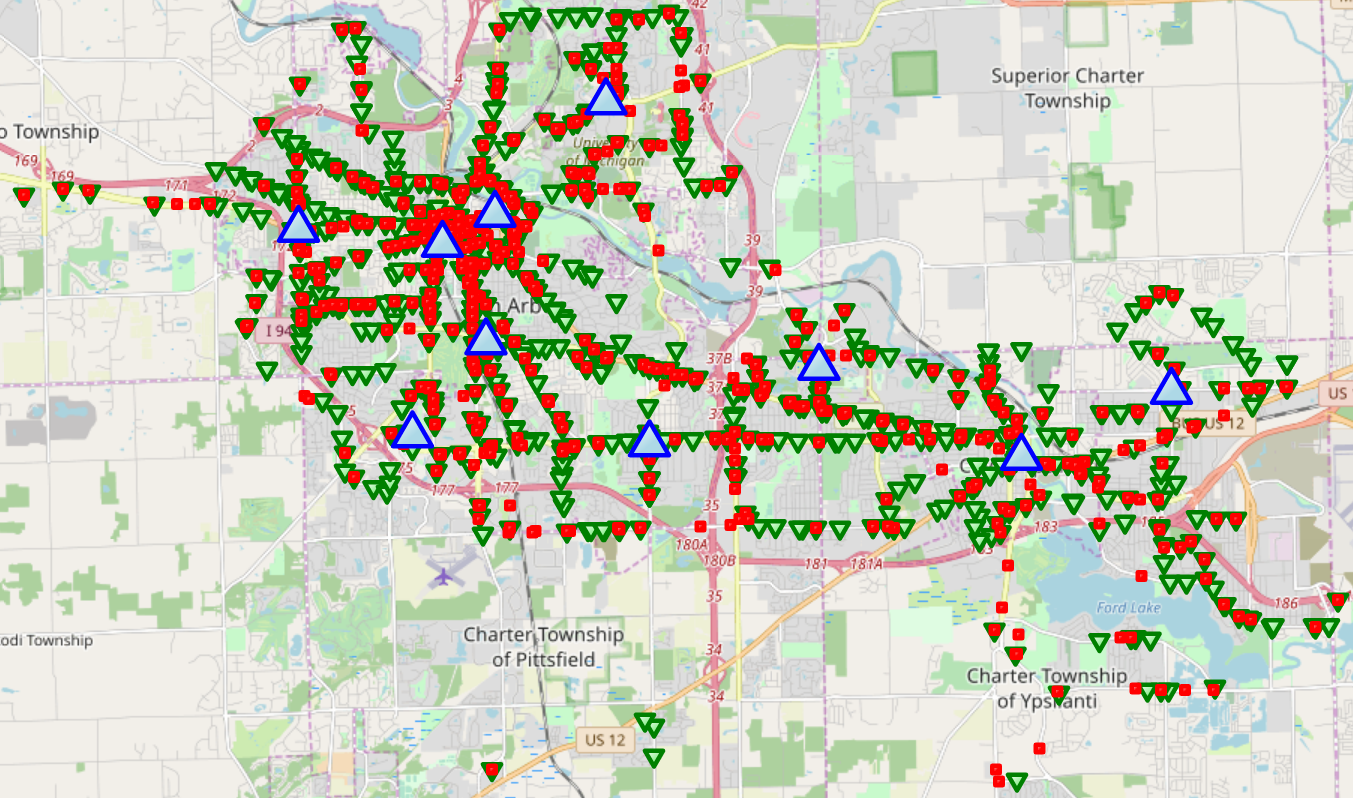}
    \caption{Illustration of the Real Case Study Instance from 6:00 am to 10:00 am: Origins (green inverted triangles), Destinations (red squares), and Hub locations (blue triangles).}
    \label{fig:instance}
\end{figure}

This section illustrates the potential of ridesharing using ridership data from 6:00am to 10:00am. Figure \ref{fig:instance} specifies the potential hub locations, and the origins and destinations of the considered riders, represented by up-oriented blue triangles, down-oriented green triangles, and small red squares respectively.

\begin{figure}[!t]
    \begin{subfigure}[b]{.5\linewidth}
        \centering
        \includegraphics[width=0.97\textwidth]{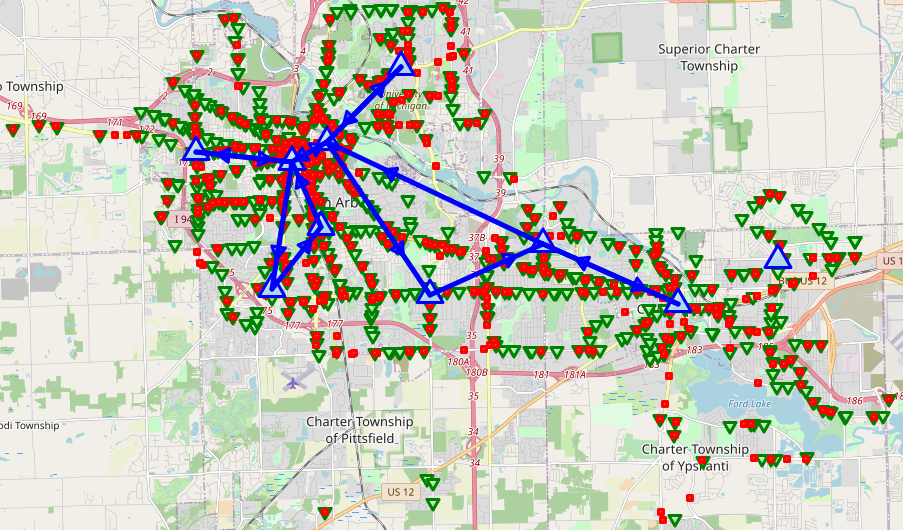}
        \caption{$K=1$}\label{fig:cap1-2}
    \end{subfigure} %
    \begin{subfigure}[b]{.5\linewidth}
        \centering
        \includegraphics[width=0.97\textwidth]{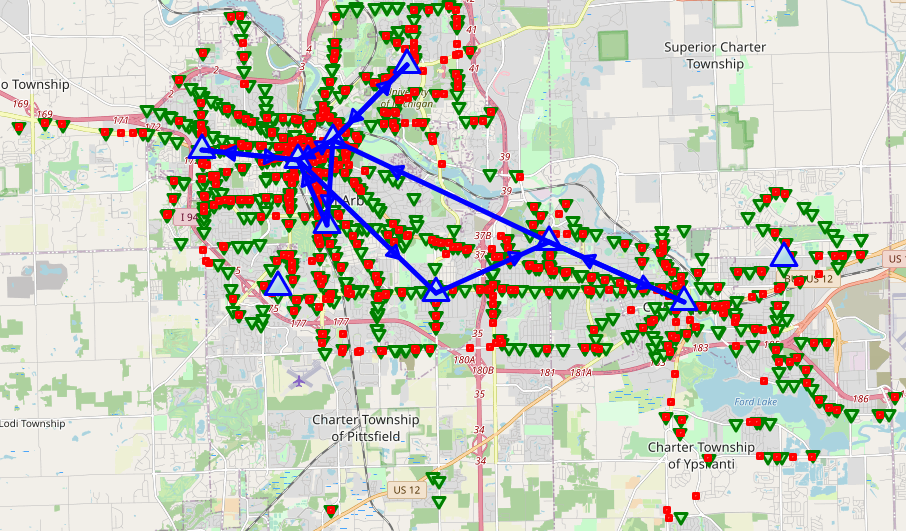}
        \caption{$K\in\{2,3,4\}$}\label{fig:cap3-4}
    \end{subfigure}
    \caption{Visualization of Optimal Design for Different Shuttle Capacities.}
    \label{fig:network_sol}
\end{figure}

Figure \ref{fig:network_sol} shows the optimal hub network obtained by solving the ODMTS design model for different shuttle capacities. Note that some of the arcs are bidirectional, while others follow a single direction. Intuitively, this is related to the spatial distribution of destinations of each commodity with respect to its origin, as well as to the weak connectivity conditions imposed by Constraints \eqref{topology}. For $K = 1$, the resulting network consists of 14 opened lines and has a large 3-hub cycle at its center that connects the two most populated areas in the extremes of the map, each extreme having in turn its own sub-network with short-length cyclical routes. When $K\in\{2,3,4\}$, however, the number of opened bus lines decreases to 13 by disconnecting one hub in the western side and resizing the associated loop. The resulting central sub-network now includes 4 bus lines that describe a 4-hub cycle that connects to both extremes of the territory. Observe that increasing the shuttle capacity results in a few modifications to the optimal network. The economies of scale of ridesharing allow shuttles to drive riders to/from hubs that are further away from their origins/destinations for a substantially lower distance cost. As a result, some bus lines that are opened when $K=1$, can be closed to achieve additional savings.

\begin{figure}[!t]
    \centering
    \includegraphics[width=0.8\textwidth]{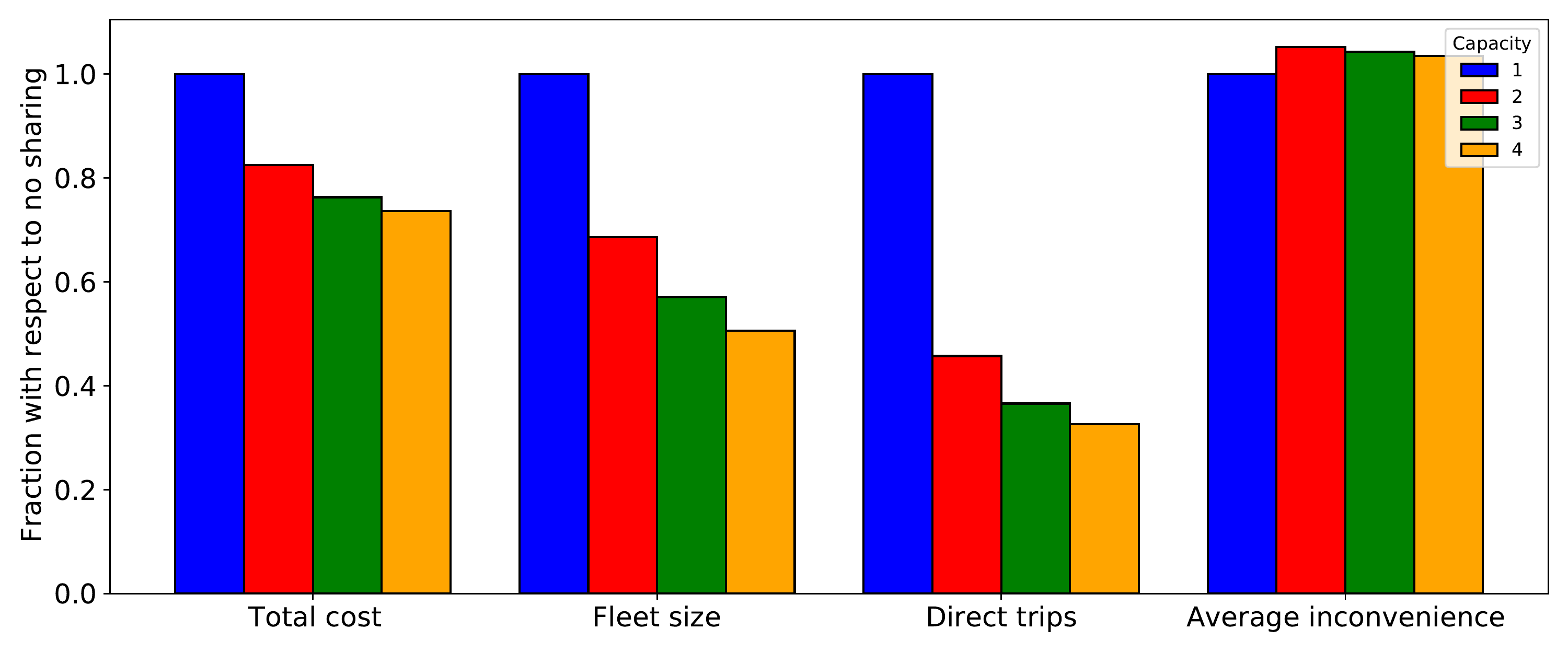}
    \caption{Relative Effect of Ridesharing on Total Cost, Number of Shuttles, Number of Direct Routes, and Average Rider Inconvenience.}
    \label{plot:total}
\end{figure}

\begin{table}[!t]
\caption{Effect of Ridesharing on Total Cost, Number of Shuttles, Frequency of Direct O-D Routes, and Number of Opened Bus Legs.}
\label{tab:total}
\centering
\begin{tabular}{ccccc}
\hline
\multirow{2}{*}{$K$} & Total & Fleet & \# direct & \# opened \\ 
& cost & size & O-D routes & bus legs \\ \hline
    1 & \$ 25,951.03 & 331 & 2,800 & 14 \\
2 & \$ 21,396.08 & 229 & 1,280 & 13  \\
3 & \$ 19,798.72 & 186 & 1,024 & 13  \\
4 & \$ 19,134.08 & 165 & 902 & 13  \\ \hline
\end{tabular}
\end{table}

\begin{table}[!t]
\caption{Effect of Ridesharing on Average Convenience (Measured in Minutes per Rider) and Average Shuttle Usage (Riders per Shuttle Ride).}
\label{tab:usage-conv}
\centering
\begin{tabular}{cccc}
\hline
\multirow{2}{*}{$K$} & Average & Average shuttle \\ 
& inconvenience & usage \\ \hline
1 & 15.31 & 1.00 \\
2 & 16.11 & 1.58 \\
3 & 15.97 & 1.89 \\
4 & 15.88 & 2.09 \\ \hline
\end{tabular}
\end{table}

\begin{figure}[!t]
    \centering
    \includegraphics[width=0.8\textwidth]{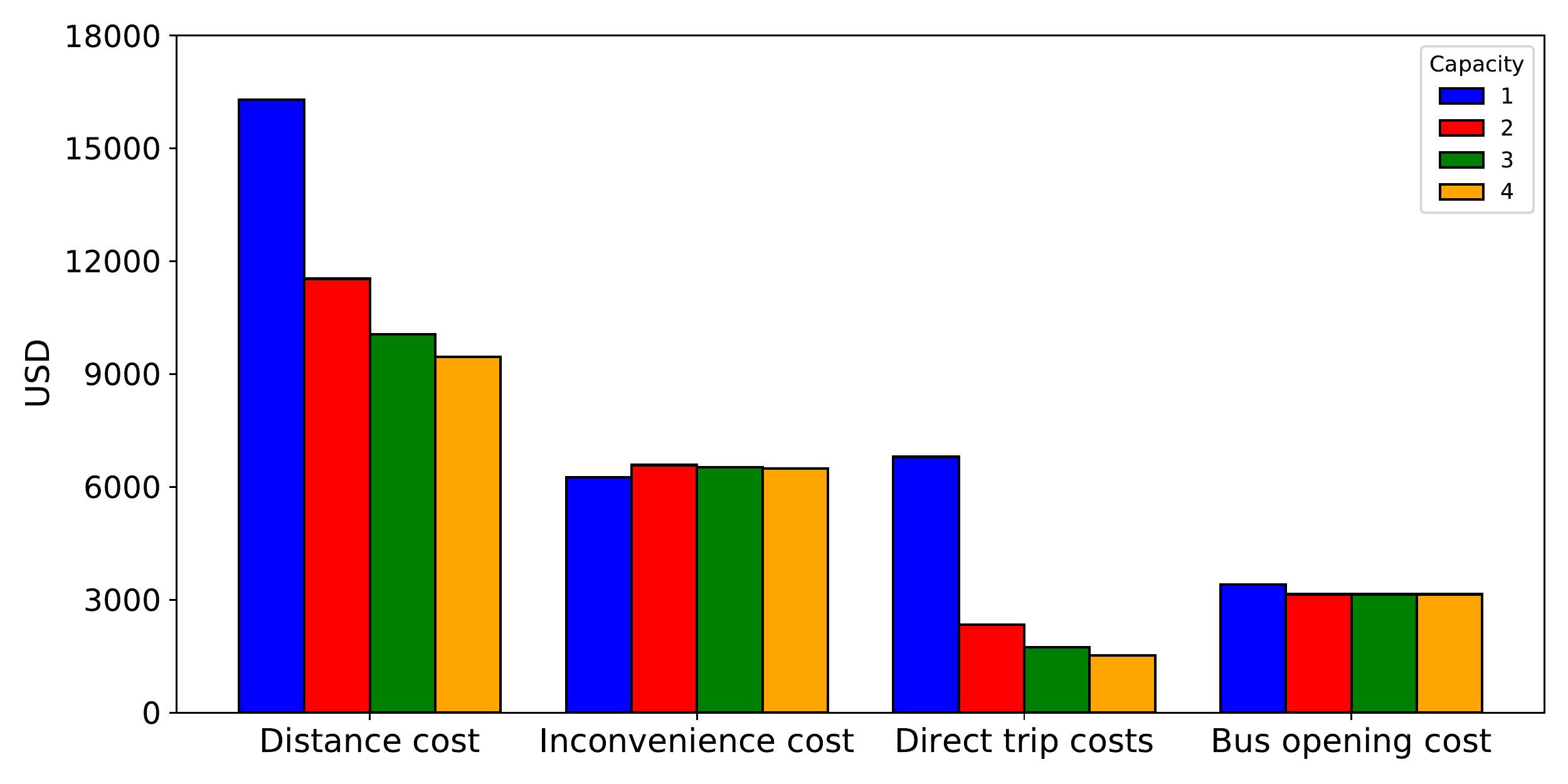}
    \caption{Effect of Ridesharing on each Type of Cost.}
    \label{plot:base-breakdown}
\end{figure}

Table \ref{tab:total} shows the effect of shuttle capacity on the total costs and the number of direct shuttle routes, and Figure \ref{plot:total} illustrates the relative variations of the associated total costs, the number of direct O-D routes, and the average inconvenience. Table \ref{tab:usage-conv} contains the average inconvenience and average shuttle usage for all the shuttle capacity values considered, and Figure \ref{plot:base-breakdown} displays a cost breakdown for the different types of costs for all capacities.

\par The results show that the total cost is reduced by 17.6\% when $K$ is increased just from 1 to 2, and by up to 26.3\% when it is further increased to 4. On the other hand, a larger shuttle capacity induces an increase in inconvenience, as observed in Table \ref{tab:usage-conv} and Figure \ref{plot:total}. However, this decrease in convenience is low, with only a 5.2\% degradation when $K$ is increased to 2. More interestingly, further increasing $K$ improves the average inconvenience, narrowing the relative degradation down to only 3.7\% when $K=4$. Unsurprisingly, a shuttle capacity of $K=1$ improves convenience; yet when $K$ is large enough, namely $K\in\{3,4\}$, it becomes beneficial to group riders in longer shared shuttle routes that drop them off (pick them up) at a hub closer to their destination (origin), saving them a number of intermediate transfers that they would incur if $K=2$. Despite these results, the average shuttle occupation is small compared to the maximum capacity $K$, being near 50\% of the shuttle capacity when $K=4$ as shown in Table \ref{tab:usage-conv}.

A similar decrease is observed in Figure \ref{plot:base-breakdown} for the costs incurred by direct O-D routes.  As $K$ increases, a major decrease in the number of direct rides is observed, going down from 2,800 when $K=1$ to only 902 when $K=4$, which constitutes a 67.8\% reduction. This in turn dramatically decreases the cost associated with direct rides, producing a 74.9\% reduction for $K=4$ compared to the cost incurred when $K=1$.

Since multiple passengers may complete their shuttle legs in a common route when ridesharing is allowed, a reduction of the number of shuttles is expected as the shuttle capacity becomes larger. Figure \ref{plot:total} and Table \ref{tab:total} present the effect of shuttle capacity on the optimal fleet size. For $K=2$, the total number of shuttles required to serve all the routes experiences a considerable decrease of 30.8\%, and these savings increase to 50.1\% when $K=4$. This illustrates the significant potential savings from adopting ridesharing since the capital expenditures for shuttle fleet can be divided by 2 when increasing the shuttle capacity. In addition, a fleet-size reduction is beneficial from a logistic, managerial, and environmental point of view, as a smaller fleet produces less traffic congestion and emission, and less coordination at pickup and dropoff locations (e.g., at bus hubs), is easier to coordinate.

\subsection{Time Bucket Sensitivity} \label{sec:exp:W}

This section analyzes the impact of the time bucket size $W$. The experiments replicate the simulation from Section \ref{sec:exp:base_case} with values $W\in\{1, 5\}$. Obviously this parameter has no effect on the results if $K=1$. Table \ref{tab:total-W} shows that decreasing $W$ to 1 minute results in a total cost increase of up to 6.2\%, whereas increasing $W$ up to 5 minutes yields a cost reduction of up to 1.8\%.  Likewise, the fleet size seems to be robust to changes in the value of $W$: decreasing $W$ to 1 minute produces a 5.5\% increase in the number of shuttles, while raising $W$ to 5 minutes results in an average decrease of 2.1\%. This is also reflected in the number of direct O-D routes: a value of $W=5$ results in a 3.8\% reduction of direct O-D routes, while $W=1$ produces an average increase of 20.4\%. The only exception to the observed pattern is the case $K=2$, where increasing $W$ from 3 to 5 minutes results in a slim 1.3\% increase in the fleet size. This is reasonable since the fleet size is not optimized by the ODMTS design model, and such a minor change may occur due to the selection of other cost-effective routes when changing the value of $W$.

\begin{table}[!t]
\caption{Effect of Ridesharing on Total Cost, Number of Shuttles, Frequency of Direct O-D Routes, and Number of Opened Bus Legs.}
\label{tab:total-W}
\centering
\begin{tabular}{cccccc}
\hline
\multirow{2}{*}{$W$} & \multirow{2}{*}{$K$} & Total & Fleet & \# direct & \# opened \\ 
& & cost & size & O-D routes & bus legs \\ \hline
\multirow{3}{*}{1} 
& 2 & \$ 22,147.80 & 235 & 1,461 & 13  \\
& 3 & \$ 20,881.17 & 195 & 1,236 & 13  \\
& 4 & \$ 20,321.92 & 180 & 1,139 & 13  \\ \hline
\multirow{3}{*}{3}
& 2 & \$ 21,396.08 & 229 & 1,280 & 13  \\
& 3 & \$ 19,798.72 & 186 & 1,024 & 13  \\
& 4 & \$ 19,134.08 & 165 & 902 & 13  \\ \hline
\multirow{3}{*}{5} 
& 2 & \$ 21,241.85 & 224 & 1,255 & 13  \\
& 3 & \$ 19,559.67 & 185 & 973 & 13  \\
& 4 & \$ 18,790.22 & 159 & 862 & 13  \\ \hline
\end{tabular}
\end{table}

Results on passenger inconvenience and average shuttle utilization are summarized in Table \ref{tab:usage-conv-W}. All changes in inconvenience due to perturbing $W$ are negligible with respect to the base case $W=3$. In general, a larger value of $W$ translates into greater shuttle utilization and fewer direct routes, which slightly increase the overall inconvenience. An exception is the case $K=4$, where the value $W=5$ is large enough so that the larger set of riders that can be grouped results in shuttle routes that are efficient in both cost and duration. In terms of shuttle utilization, decreasing $W$ to 1 minute reduces the average number of riders in a route by 7.7\%, whereas increasing $W$ to 5 minutes results in an overall increase of 1.9\%.

\begin{table}[!t]
\caption{Effect of Ridesharing on Average Inconvenience (Measured in Minutes per Passenger) and Average Shuttle Usage (Passengers per Shuttle Ride).}
\label{tab:usage-conv-W}
\centering
\begin{tabular}{cccc}
\hline
\multirow{2}{*}{$W$} & \multirow{2}{*}{$K$} & Average & Average shuttle \\ 
& & inconvenience & usage \\ \hline
\multirow{3}{*}{1}
& 2 & 15.89 & 1.50 \\
& 3 & 15.83 & 1.75 \\
& 4 & 15.67 & 1.87 \\ \hline
\multirow{3}{*}{3}
& 2 & 16.11 & 1.58 \\
& 3 & 15.97 & 1.89 \\
& 4 & 15.88 & 2.09 \\ \hline
\multirow{3}{*}{5}
& 2 & 16.11 & 1.59 \\
& 3 & 15.99 & 1.94 \\
& 4 & 15.77 & 2.14 \\ \hline
\end{tabular}
\end{table}

For each shuttle capacity value $K$, the considered values $W\in\{1,5\}$ results in an optimal bus network which is identical to the one obtained for $W=3$ in Section \ref{sec:exp:base_case}. This evidences the robustness of the bus network design with respect to both the shuttle capacity and the length of the time buckets in which multiple riders can be grouped in a single shuttle route.

\subsection{Sensitivity Analysis on the Shuttle Route Duration Threshold} \label{sec:delta_sensitivity}

This section assesses the impact of the route duration threshold $\delta$. Results for the different metrics are summarized in Tables \ref{tab:delta_sensitivity_1} and \ref{tab:delta_sensitivity_2}, assuming a capacity $K=3$ and a time bucket length $W=3$. The results show that, even after considerably increasing the threshold up to 150\%, the total cost only improves by a mere 2.5\%, while the average inconvenience increases by 1.3\%. As expected, increasing the threshold creates additional opportunities of grouping people in shared routes, decreasing the required fleet size and the number of direct O-D routes, and increasing the average shuttle occupancy. These changes do not affect the optimal network design with respect to the topology obtained for $\delta = 50\%$.

The low magnitude of these changes can be explained by the assumption that $H_r^-$ and $H_r^+$ only comprise the 3 closest hubs to the origin and destination of each commodity $r\in C$, respectively, thus not offering much more possibility of consolidation than the ones already possible with a value of $\delta = 50\%$. Increasing $\delta$ may offer further benefits than the ones currently observed for larger sets $H_r^-$ and $H_r^+$, however this would come at the cost of increased computational challenges for the MIP solver.
\begin{table}[!t]
    \centering
    \begin{tabular}{ccccc}
         \multirow{2}{*}{$\delta$} & Total & \# opened & Fleet & \# direct\\
         & cost & bus lines & size & shuttle routes \\\hline
         50\% & \$19,798.72 &  13 & 186 & 1,024
         \\
         100\% & \$19,419.28 & 13 & 172 & 941
         \\ 
         150\% & \$19,303.88 & 13 & 170 & 928
         \\ \hline
    \end{tabular}
    \caption{Effect of Shuttle Route Duration Threshold on Total Cost, Number of Shuttles, Frequency of Direct O-D Routes, and Number of Opened Bus Legs, for $(K, W) = (3, 3)$.}
    \label{tab:delta_sensitivity_1}
\end{table}

\begin{table}[!t]
    \centering
    \begin{tabular}{ccc}
         \multirow{2}{*}{$\delta$} & Average & Average \\
         & inconvenience (min) & shuttle usage \\\hline
         50\% & 15.97 & 1.89 
         \\
         100\% & 16.15 & 2.03 
         \\ 
         150\% & 16.18 & 2.08 
         \\ \hline
    \end{tabular}
    \caption{Effect of Shuttle Route Duration Threshold on Average Inconvenience (Measured in Minutes per Rider) and Average Shuttle Usage (Passengers per Shuttle Ride), for $(K, W) = (3, 3)$.}
    \label{tab:delta_sensitivity_2}
\end{table}

\subsection{Sensitivity Analysis on the Estimated Hub Arrival Time} \label{sec:exp:t0+}

\begin{figure}[!t]
    \centering
    \includegraphics[width=0.5\textwidth]{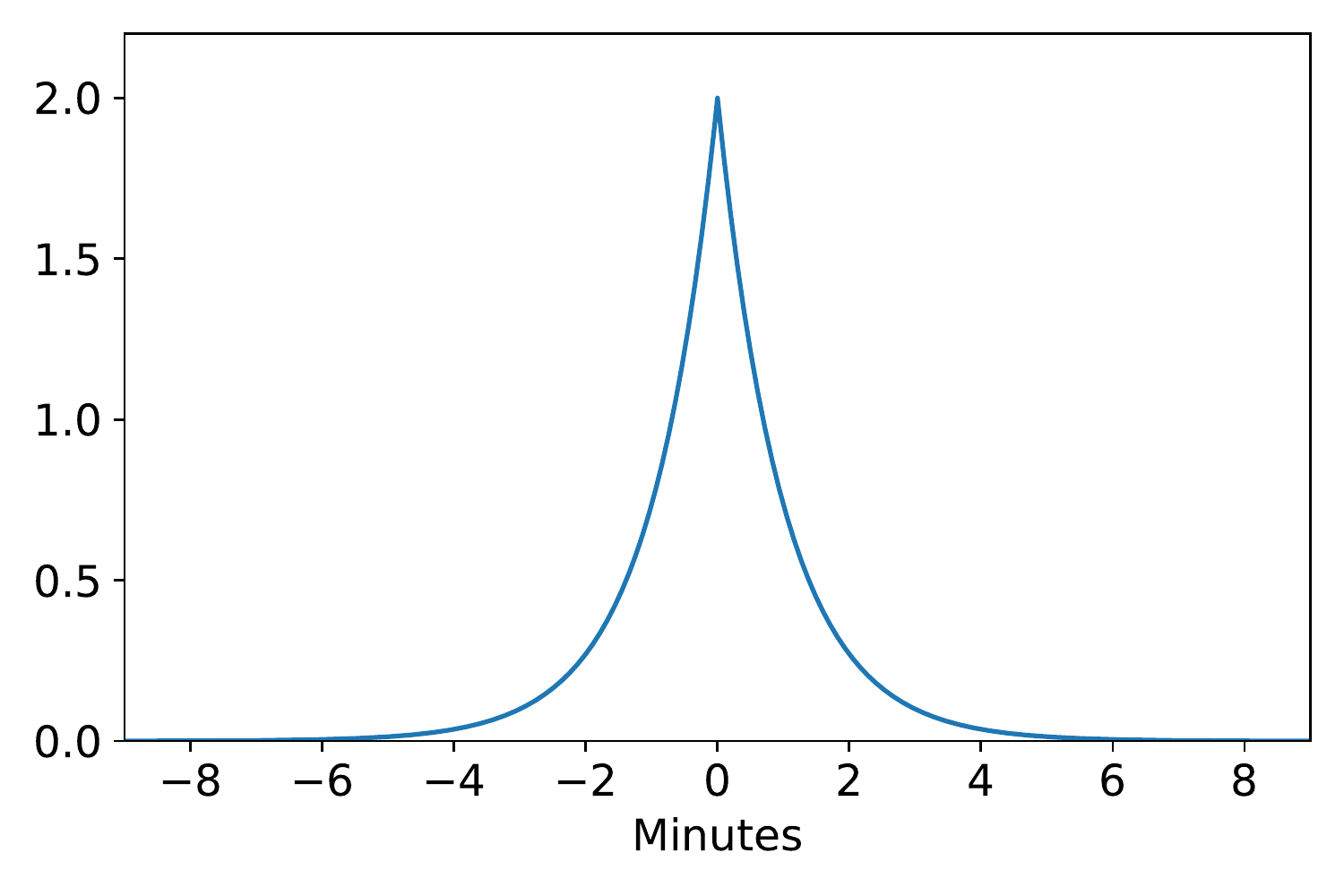}
    \caption{Density function used to sample the perturbations for $t_1(r, h)$, which corresponds to $Laplace(0, 1)$.}
    \label{fig:laplace-dist}
\end{figure}

This section studies how sensitive the proposed model is to perturbations in the estimate arrival time to the last hub $t_1(r,h)$. This analysis helps assessing the validity of using this estimation as an input instead of leaving $t_1(r,h)$ as part of the variables, which would make the model much harder. For this purpose, for each commodity $r\in C$ and each hub $h$, a noise sampled from a Laplace distribution (in minutes) is added to $t_1(r, h)$ (see Figure \ref{fig:laplace-dist} for the exact distribution).  The ODMTS design model and the sparse fleet-sizing model are then solved using the perturbed estimates. Such change in the arrival times to the last hub will result in some passengers arriving earlier or later than in the base instance from Section \ref{sec:exp:base_case}, possibly modifying the set of trips that can be consolidated in the last shuttle leg. In order to capture the effect of variations in $t_1(r,h)$, such procedure is repeated a total of 50 times and report some statistics for various metrics for shuttle capacity of $K=3$.

The results are shown in Table \ref{tab:t0+sens_res}, where the performance metrics for the perturbed instances are compared with those of the base instance. Overall, the metric values from the base instance are either contained in, or very close to, the reported range from the perturbed instances. In particular, the model proves to be robust in terms of operational cost, with a minor increase between 1.2\% to 1.6\% with respect to the base cost. Furthermore, it is also robust in terms of the inconvenience and optimal fleet size: perturbed inconvenience experiences an overall increase between -0.9\% to 1.8\% from the base inconvenience, and the perturbed optimal fleet size between -1.1\% and 3.8\%. In terms of shuttle occupancy, perturbing $t_1(r,h)$ produces an overall decrease of 8.3\% in last leg routes: the perturbations restrict the consolidation opportunities in the last leg of trips, in turn increasing the overall costs due to having a larger driven distance. This also makes long last shuttle legs too costly since the driving cost is split among fewer riders, in turn requiring riders travel more by bus; as a result, some instances show an overall increase in inconvenience. A slight overall increase of the number of direct routes is observed, which explains the raise in total cost. In the particular case of the optimal bus network, the 50 runs open exactly the same bus lines, giving additional evidence of the robustness of the model to changes in $t_1(r, h)$ and validating the assumption on its estimation.

\begin{table}[!t]
  \caption{Sensitivity Analysis of all Metrics on $t_1(r, h)$.}
\label{tab:t0+sens_res}
\centering
\begin{tabular}{lcccc}
& \multicolumn{3}{c}{Perturbed instances} & \\ \cmidrule{2-4}
\cmidrule(l){2-3}
& Min & Avg & Max & Base \\ \hline
Total cost (\$) & 20,055.14 & 20,094.82 & 20,134.56 & 19,798.72 \\
\# opened bus lines & 13 & 13 & 13 & 13 \\
Inconvenience (minutes) & 15.82 & 16.03 & 16.27 & 15.97 \\
Average shuttle usage (first leg) & 2.16 & 2.17 & 2.18 & 2.17 \\
Average shuttle usage (last leg) & 1.53 & 1.55 & 1.60 & 1.68 \\
\# direct O-D routes & 1,039 & 1,079.6 & 1,104 & 1,024 \\ 
Fleet size & 184 & 188.4 & 193 & 186 \\\hline
\end{tabular}
\end{table}

\subsection{Benefits of the ODMTS}
\label{sec:exp:transit_system_eval}

This section compares the system designed by the ODMTS design Model \eqref{mip_model}, referred to as the \textit{proposed system}, with the current public transit system. The comparison is performed under the following considerations:
\begin{enumerate}[(i)]
    \item both systems are evaluated for 6:00am--10:00pm;
    \item the proposed system is designed with the same parameter values used in Section \ref{sec:exp:base_case} and assuming a fleet of shuttles of capacity $K=3$;
    \item the comparison between the current and proposed systems uses two metrics: the daily operational costs and the average inconvenience of the passengers requesting service;
    \item the proposed system is constructed by independently solving the ODMTS design model for each of the four 4-hour time periods between 6:00am and 10:00 pm. The system total operational cost is then computed as the sum over the four obtained solutions; similarly, the average inconvenience is computed as the ratio between the combined travel time over all the riders that requested service in any of the considered time periods, and the total number of such riders.
\end{enumerate}

\begin{table}[!t]
\caption{Comparison between current and proposed public transit systems}
\label{tab:comparison-systems}
\centering
\begin{tabular}{ccc}
\hline
\multirow{2}{*}{System} & Daily total & Average  \\
& cost (\$) & inconvenience (min) \\ \hline
Current  & 105,865.85 & 24.10 \\
Proposed $(K=3)$ & 68,745.52 & 14.86 \\\hline
\end{tabular}
\end{table}

\begin{figure}[!t]
    \begin{subfigure}[b]{.5\linewidth}
        \centering
        \includegraphics[width=0.97\textwidth]{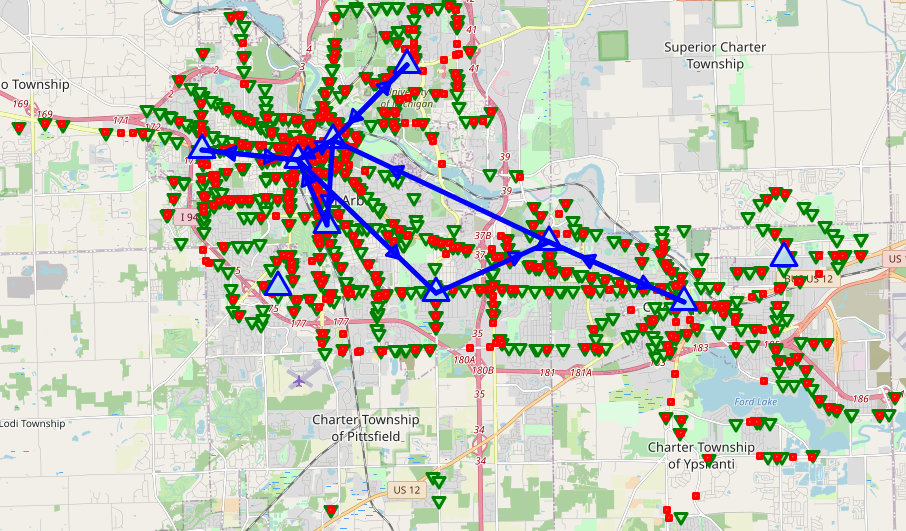}
        \caption{6:00 am - 10:00 am}\label{fig:map6}
    \end{subfigure} %
    \begin{subfigure}[b]{.5\linewidth}
        \centering
        \includegraphics[width=0.97\textwidth]{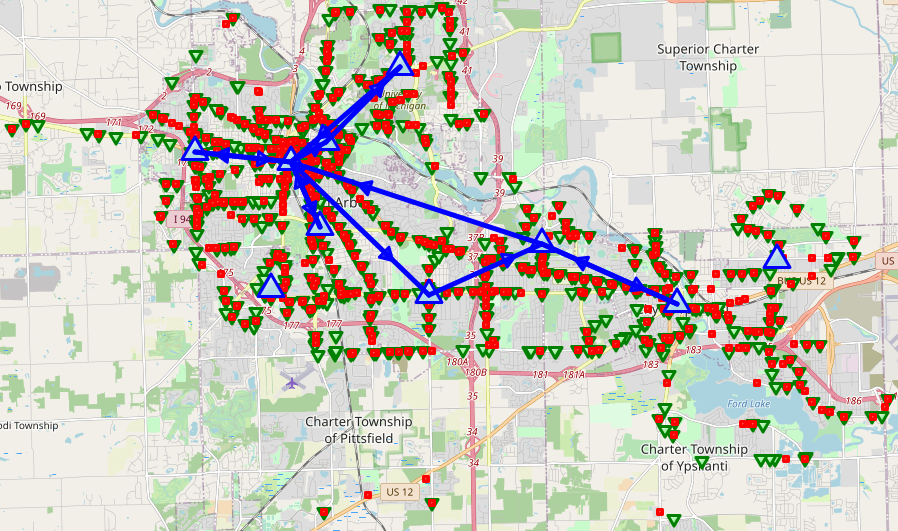}
        \caption{10:00 am - 2:00 pm}\label{fig:map10}
    \end{subfigure}

    \begin{subfigure}[b]{.5\linewidth}
        \centering
        \includegraphics[width=0.97\textwidth]{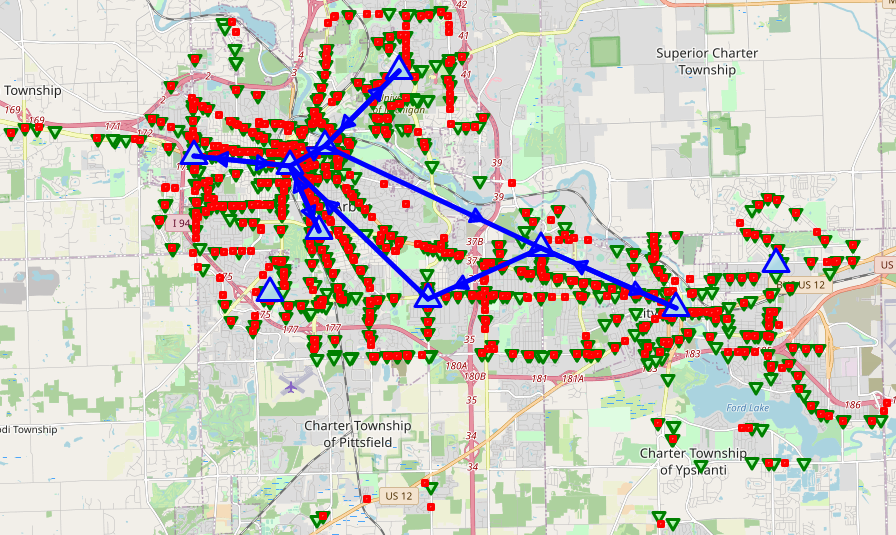}
        \caption{2:00 pm - 6:00 pm}\label{fig:map14}
    \end{subfigure} %
    \begin{subfigure}[b]{.5\linewidth}
        \centering
        \includegraphics[width=0.97\textwidth]{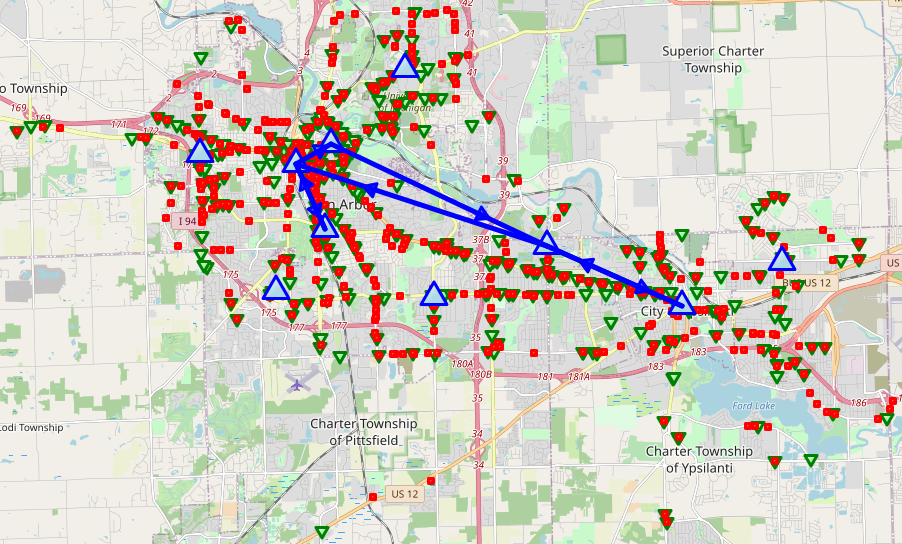}
        \caption{6:00 pm - 10:00 pm}\label{fig:map18}
    \end{subfigure}
    \caption{Visualization of Optimal Bus Networks at Different Times of an Operating Day. Figure created with \cite{OpenStreetMap}}
    \label{fig:maps-day}
\end{figure}

Table \ref{tab:comparison-systems} compares both systems in terms of operational cost and inconvenience. The results show that the proposed system achieves a total daily operational savings of \$37,120.33 a day. Although the proposed ODMTS incurs an extra cost of operating shuttles, the cost savings from using fewer buses reduce the total costs by 35.1\%. Moreover, the proposed system greatly improves the average inconvenience, reducing it by 38.3\%. Both improvements are explained by the significant reduction in the number of buses, the shuttle mode for first/last miles, and ridesharing. These findings demonstrate the great potential of an ODMTS when ridesharing is allowed.

Figure \ref{fig:maps-day} displays the optimal networks for each of the 4-hour time blocks, Note that the network barely changes from 6:00am to 6:00pm (Figures \ref{fig:map6} - \ref{fig:map14}): the only modifications are a few bus lines opened in the Western side of the territory, and the reduction of the 4-hub loop at the center to a smaller 3-hub loop between 10:00am and 2:00pm. This is in contrast to the optimal network design in the off-peak period between 6:00pm and 10:00pm where most hubs become disconnected, due to a significant decrease in demand as shown in Figure \ref{fig:map18}. Despite the reduced demand in the last 4-hour period, the weak connectivity constraints are sufficient to obtain a set of connected bus lines throughout the whole day.

\section{Conclusion} \label{sec:concl}

This work studied how to integrate ridesharing in the design of ODMTS, i.e., the possibility of grouping multiple riders in the shuttles that serve the first and last miles. The paper addressed two gaps in existing tools for designing ODMTS. First, the paper included ridesharing in the shuttle rides. Second, it proposed novel fleet-sizing algorithms for determining the number of shuttles needed to meet the performance metrics of the ODMTS design. Both contributions were based on MIP models. For the ODMTS design, the MIP featured decision variables representing whether to select a pickup or dropoff route grouping riders who travel to/from the same hub. The fleet-sizing optimization was modeled as a minimum flow problem with covering constraints. The natural formulation leads to a dense graph and significant computational issues, while the reformulation leads to a sparse graph.

The proposed framework was applied to a real case study involving a public transit system of the broader Ann Arbor and Ypsilanti region in Michigan, USA. An extensive computational analysis measured the impact of ridesharing on the design of an ODMTS using data from the busiest 4-hour time period on a weekday. It was observed that ridesharing can reduce costs by more than 25\% when using shuttles with capacity 4 (compared to shuttles with capacity 1), in exchange of a slight increase of around 4\% in transit times. Additionally, the study shows that ridesharing yields a considerable reduction of the minimum fleet size required to serve all shuttle routes: the reductions range from 29.2\% when using shuttles of capacity 2 to 45.3\% for shuttles of capacity 4. Interestingly, shuttles are rarely used at full capacity with an average occupancy of 2.1 passengers per route for shuttles of capacity 4.

The paper also conducted a sensitivity analysis on the time window used for consolidation and the estimation of the arrival time to the last hub of riders, showing that the results are robust overall. The proposed ODMTS was also compared to the existing public transit system in terms of cost and convenience. The findings suggest that the ODMTS reduces cost by 35\% and transit times by 38\%.

There are several interesting directions for future research. First of all, it should be noted that the numerical results shown in this paper provide an optimistic bound of the potential benefits of ridesharing in terms of cost and fleet size, as the conducted experiments involve perfect knowledge of the transportation requests. This is a reasonable assumption in general since transit riders are overall highly loyal and predictable. However, a possible research direction would consist of determining the advantage of ridesharing in settings where the demand is revealed over time. From an algorithmic perspective, our current research focuses on optimization techniques that scale better with the number of hubs $|H|$ and the shuttle route duration threshold $\delta$. Incorporating lateness and uncertainty in the fleet-sizing algorithm is an important extension, which makes the problem significantly more challenging as a single route may now be served at different start times depending on the route served immediately before. The fleet-sizing model can also consider additional objectives that depend on the task sequencing, e.g., travel distance or travel time. These objectives do not admit the pre-filtering step employed to sparsify the network and hence raise interesting computational issues. From a practical point of view, the integration of mode choice models (e.g., \cite{basciftci2021capturing}), additional modes of transportation (e.g., a rail system), modeling of more complex bus lines with intermediate stops, and the consideration of shuttle routes that combine pickups, dropoffs, direct trips, and repositioning may offer additional invaluable insights.

\par
\subsection*{Acknowledgements}

Ramon Auad is supported by the National Agency for Research and
Development (ANID) through the scholarship program Doctorado Becas
Chile 2017 - 72180404. The work was partly supported by NSF Leap HI
proposal NSF-1854684 and Department of Energy Research Award 7F-30154.

\bibliography{references, fleet_sizing}
\bibliographystyle{apalike}
\vfill
\pagebreak

\section*{Appendices}
\appendix
\section{Shuttle Scheduling} \label{appA}
Given an optimal solution $v^*$ obtained from solving Model \eqref{mip:fleet_size_tuned}, Algorithm \ref{recover_shuttle_paths} allows to obtain a schedule to be followed by each shuttle. Note that each path from source $s$ to sink $s'$ described by $v^*$ specifies a sequence of visited nodes (with each node corresponding to a shuttle route), which ultimately can be translated into a shuttle schedule. The key idea of Algorithm \ref{recover_shuttle_paths} is then to iteratively identify existing paths (i.e., set of connected arcs whose flow value is non-zero) connecting $s$ and $s'$ in the solution, and assign the nodes/routes covered by each path to a different shuttle. More specifically, when a path is found, a shuttle is scheduled to serve all the routes present in the path that have not been already covered by a previous shuttle, and then the flow value of each arc in the identified path is reduced by 1 unit; this last step is required since the flow variables in the sparse formulation by Model \eqref{mip:fleet_size_tuned} are unconstrained, and not doing it would end up assigning a route to multiple shuttles. At termination, Algorithm \ref{recover_shuttle_paths} effectively assigns a valid schedule of routes to each shuttle, as it finds a number of non-empty schedules that exactly matches the optimal number of shuttles.

\begin{algorithm}[H]
  \caption{Recovery of Shuttle Schedules}
  \label{recover_shuttle_paths}
\begin{algorithmic}[1]
\Require Optimal flow vector $v^*$ from the sparse fleet-sizing model.
\Ensure Set of shuttle paths $P$.
\State $P\gets \emptyset$
\State $\overline{\Omega} \gets \Omega$
\While{$v^* \neq 0$}
    \State $p\gets\emptyset$
    \State Identify a path $\mathcal{A}_{path} = \{(s,\omega_1), (\omega_1,\omega_2), \dots, (\omega_{k}, s')\}$ such that $v^*_{a} > 0, \forall a\in\mathcal{A}_{path}$.
    \For{$a\in\mathcal{A}_{path}$}
        \State $v^*_a \gets v^*_a - 1$
    \EndFor
    \For{$j\in\{1,2,\dots,k\}$}
        \If{$\omega_j \in \overline{\Omega}$}
            \State $p\gets p\cup\{\omega_j\}$
            \State $\overline{\Omega} \gets \overline{\Omega} \setminus \{\omega_j\}$
        \EndIf
    \EndFor
    \State $P\gets P\cup\{p\}$
\EndWhile
\end{algorithmic}
\end{algorithm}

\end{document}